\documentclass[manyauthors]{fundam}

\usepackage{url, amssymb, enumerate, epsfig, subcaption, mathrsfs} 
\usepackage[ruled,lined]{algorithm2e}
\usepackage{graphicx}
\usepackage{tikz}
\usetikzlibrary{automata, positioning, arrows}

\usepackage{todonotes}

\usepackage{hyperref}

\begin{document}


\setcounter{page}{197}
\publyear{24}
\papernumber{2179}
\volume{191}
\issue{3-4}

\finalVersionForARXIV


\title{On Three Domination-based Identification Problems \\  in Block Graphs}

\author{Dipayan Chakraborty\thanks{Also works: Department of Mathematics and Applied Mathematics,
                              University of Johannesburg, Auckland Park, 2006, South Africa}\thanks{Address for correspondence:
                                       LIMOS, 1 rue de la Chebarde, Campus des C\'ezeaux, 63178 Aubi\`ere Cedex, France}
\\
Universit\'e Clermont-Auvergne, CNRS\\
Mines de Saint-\'Etienne, Clermont-Auvergne-INP\\
 LIMOS, 63000 Clermont-Ferrand, France \\
dipayan.chakraborty{@}uca.fr
\and
Florent Foucaud\\
Universit\'e Clermont-Auvergne, CNRS\\
 Mines de Saint-\'Etienne, Clermont-Auvergne-INP\\
 LIMOS, 63000 Clermont-Ferrand, France \\
florent.foucaud{@}uca.fr
\and
Aline Parreau
\\
Univ Lyon, CNRS, INSA Lyon, UCBL\\
 Centrale Lyon, Univ Lyon 2, LIRIS, UMR5205\\
 69622 Villeurbanne Cedex, France \\
aline.parreau{@}univ-lyon1.fr
\and
Annegret Wagler\\
Universit\'e Clermont-Auvergne, CNRS\\
 Mines de Saint-\'Etienne, Clermont-Auvergne-INP\\
  LIMOS, 63000 Clermont-Ferrand, France \\
annegret.wagler{@}uca.fr
}

\maketitle

\runninghead{D. Chakraborty et al.}{On Three Domination-based Identification Problems in Block Graphs}

\vspace*{-4mm}
\begin{abstract}
  The problems of determining the minimum-sized \emph{identifying}, \emph{locating-dominating} and \emph{open locating-dominating codes} of an input graph are special search problems that are challenging from both theoretical and computational viewpoints. In these problems, one selects a dominating set $C$ of a graph $G$ such that the vertices of a chosen subset of $V(G)$ (i.e. either $V(G)\setminus C$ or $V(G)$ itself) are uniquely determined by their neighborhoods in $C$. A typical line of attack for these problems is to determine tight bounds for the minimum codes in various graph classes. In this work, we present tight lower and upper bounds for all three types of codes for \emph{block graphs} (i.e. diamond-free chordal graphs). Our bounds are in terms of the number of maximal cliques (or \emph{blocks}) of a block graph and the order of the graph. Two of our upper bounds verify conjectures from the literature — with one of them being now proven for block graphs in this article. As for the lower bounds, we prove them to be linear in terms of both the number of blocks and the order of the block graph. We provide examples of families of block graphs whose minimum codes attain these bounds, thus showing each bound to be tight.

\end{abstract} 

\section{Introduction}

For a graph $G$ that models a facility or a multiprocessor network, detection devices can be placed at its vertices to locate an intruder (like a faulty processor, a fire or a thief).
Depending on the features of the detection devices, different types of dominating sets can be used to determine the optimum distributions of these devices across the vertices of $G$. In this article, we study three problems arising in this context, namely the three types of dominating sets — called the \emph{identifying codes}, \emph{locating-dominating codes} and \emph{open locating-dominating codes} — of a given graph. Each of these problems has been extensively studied during the last decades (see the bibliography maintained by Lobstein~\cite{lobstein2012watching}). These three types of codes are among the most prominent notions within the larger research area of identification problems in discrete structures pioneered by R\'enyi~\cite{renyi1961}, with numerous applications, for example in fault-diagnosis~\cite{Rao93}, biological testing~\cite{MS85} or machine learning~\cite{CN98}.

\begin{figure}[t!]
     \centering
         \includegraphics[scale=0.9]{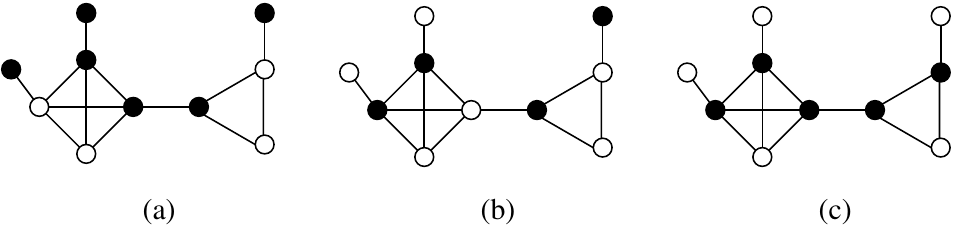}\vspace*{-1mm}
     \caption{Examples of (a) an ID-code, (b) an LD-code and (c) an OLD-code. The set of black vertices in each of the three graphs constitute the respective code of the graph.}     \label{fig_ID-OLD-LD Examples}\vspace*{-2mm}
\end{figure}

\medskip
Let $G=(V(G),E(G))$ be a graph, where $V(G)$ and $E(G)$ denote the set of vertices (also called the \emph{vertex set}) and the set of edges (also called the \emph{edge set}), respectively, of $G$. The \emph{(open) neighborhood} of a vertex $u \in V(G)$ is the set $N_G(u)$ of all vertices of $G$ adjacent to $u$; and the set $N_G[u] = \{u\} \cup N_G(u)$ is called the \emph{closed neighborhood} of $u$. A subset $C \subseteq V(G)$ is called an \emph{identifying code} \cite{karpovsky1998new} (or an \emph{ID-code} for short) of $G$ if %
\begin{itemize}
	\item $N_G[u]\cap C\neq \emptyset$ for all vertices $u \in V(G)$ (i.e. $C$ is said to be a \emph{dominating set} of $G$, or is said to possess the property of \emph{domination} in $G$); and
	
	\item $N_G[u]\cap C\neq N_G[v]\cap C$ for all distinct vertices $u,v \in V(G)$ (i.e. $C$ is called a \emph{closed-separating set} of $G$, or is said to possess the property of \emph{closed-separation} in $G$).
\end{itemize}
See Figure~\ref{fig_ID-OLD-LD Examples}a for an example of an ID-code. A graph $G$ admits an ID-code if and only if $G$ has no \emph{closed-twins} (i.e. a pair of distinct vertices $u,v\in V$ with $N_G[u]=N_G[v]$). Said differently, a graph $G$ admits an ID-code if and only if $G$ is \emph{closed-twin-free}.

\medskip
A subset $C \subseteq V(G)$ is called a \emph{locating-dominating code} \cite{slater1987domination, slater1988dominating} (or an \emph{LD-code} for short) of $G$ if
\begin{itemize}
	\item $N_G[u]\cap C\neq \emptyset$ for all vertices $u \in V(G)$ (i.e. $C$ is a dominating set of $G$); and
	
	\item $N_G(u)\cap C\neq N_G(v)\cap C$ for all distinct vertices $u,v \in V(G) \setminus C$ (i.e. $C$ is called a \emph{locating set} of $G$, or is said to possess the property of \emph{location} in $G$).
\end{itemize}
See Figure \ref{fig_ID-OLD-LD Examples}b for an example of an LD-code. Note that every graph has an LD-code.

\medskip
Finally, a subset $C \subseteq V(G)$ is called an \emph{open locating-dominating code} \cite{seo2010open} (or an \emph{OLD-code} for short) of $G$ if
\begin{itemize}
	\item $N_G(u)\cap C\neq \emptyset$ for all vertices $u \in V(G)$ (i.e. $C$ is called an \emph{open-dominating set} of $G$, or is said to possess the property of \emph{open-domination} in $G$)\footnote{The property of open-domination is often called \emph{total-domination} in the literature. See for example \cite{henning2013total}}; and
	
\eject
	\item $N_G(u)\cap C\neq N_G(v)\cap C$ for all distinct vertices $u,v \in V(G)$ (i.e. $C$ is called an \emph{open-separating set} of $G$, or is said to possess the property of \emph{open-separation} in $G$).
\end{itemize}
See Figure \ref{fig_ID-OLD-LD Examples}c for an example of an OLD-code. A graph $G$ admits an OLD-code if and only if $G$ has neither isolated vertices nor \emph{open-twins} (i.e. a pair of distinct vertices $u,v\in V(G)$ such that $N_G(u)=N_G(v)$). Again, said differently, a graph $G$ admits an OLD-code if and only if $G$ has no isolated vertices and is \emph{open-twin-free}.

\medskip

A graph with neither open- nor closed-twins is simply referred to as \emph{twin-free}.

For the rest of this article, we often simply use the word \emph{code} to mean any of the above three ID-, LD- or OLD-codes without distinction. Given a graph $G$, the \emph{identifying code number} $\gamma^{ID}(G)$ (or \emph{ID-number} for short), the \emph{locating-dominating number} $\gamma^{LD}(G)$ (or \emph{LD-number} for short) and the \emph{open locating-dominating number} $\gamma^{OLD}(G)$ (or \emph{OLD-number} for short) of a graph $G$ are the minimum cardinalities among all ID-codes, LD-codes and OLD-codes, respectively, of $G$. In other words, for simplicity, for any symbol X $\in \{$ID, LD, OLD$\}$, we have the X-number: $\gamma^X (G) = \min \{ |C| : C \text{ is an X-code of } G \}$. In the case that all three codes are addressed together as one unit anywhere in the text, i.e. any specific symbol for X $\in \{$ID, LD, OLD$\}$ is irrelevant to the context, we then simply refer to the X-numbers as the \emph{code numbers} of $G$.

\medskip
Given two sets $A$ and $B$, the set $A \triangle B = (A \setminus B) \cup (B \setminus A)$ is called the \emph{symmetric difference} of $A$ and $B$. For a subset $C \subset V(G)$ and distinct vertices $u,v \in V(G)$, if there exists a vertex $w \in (N_G(u) \cap C) \triangle (N_G(v) \cap C)$ (resp. $(N_G[u] \cap C) \triangle (N_G[v] \cap C)$), then $w$ and $C$ are said to \emph{open-separate} (resp. \emph{closed-separate}) the vertices $u$ and $v$ in $G$.

\subsection{Known results}

Given a graph $G$, determining $\gamma^{ID}(G)$ or $\gamma^{LD}(G)$ is, in general, NP-hard \cite{charon2003minimizing} and remains so for several graph classes where other hard problems become easy to solve. These include bipartite graphs~\cite{charon2003minimizing} and two subclasses of chordal graphs, namely split graphs and interval graphs~\cite{foucaud2017identification}. In fact, for both bipartite and split graphs, it is NP-hard to even approximate the ID-number and LD-number within a factor of $\log |V(G)|$~\cite{DBLP:journals/jda/Foucaud15}. Determining $\gamma^{OLD}(G)$ is also, in general, NP-hard~\cite{seo2010open} and remains so for perfect elimination bipartite graphs~\cite{P15} and interval graphs~\cite{foucaud2017identification}. On the other hand, determining $\gamma^{OLD}(G)$ becomes APX-complete for chordal graphs with maximum degree $4$~\cite{P15}.

\medskip
As these problems are computationally very hard, a typical line of attack is to determine bounds on the code numbers for specific graph classes. Closed formulas for these parameters have so far been found only for restricted graph families (e.g. for paths and cycles~\cite{slater1988dominating, seo2010open, bertrand2004identifying}, for stars~\cite{DBLP:journals/dm/GravierM07}, for complete multipartite graphs~\cite{argiroffo2018polyhedra,argiroffo2022polyhedra} and for some subclasses of split graphs including thin headless spiders~\cite{argiroffo2014study}). Lower bounds for all three code numbers for several graph classes like interval graphs, permutation graphs, cographs~\cite{DBLP:journals/tcs/FoucaudMNPV17} and lower bounds for ID-numbers for trees~\cite{bertrand20051}, line graphs~\cite{foucaud2013identifying}, planar graphs~\cite{Rall1984location} and many others of bounded VC-dimension~\cite{bousquet2015identifying} have been determined. As far as upper bounds for the code numbers are concerned, for certain graph classes, upper bounds for ID-codes (see~\cite{BFH15,FL22,FP12}), LD-codes (see \cite{BFH15,foucaud2016location,garijo2014difference}) and OLD-codes (see \cite{HY14}) have been obtained.

\subsection{Our work}

In this paper, we consider the family of block graphs, defined by Harary in~\cite{harary1963}, see also~\cite{H79} for equivalent characterizations. A \emph{block graph} is a graph in which every maximal 2-connected subgraph is complete. In a block graph, every maximal complete subgraph is called a \emph{block}. Equivalently, block graphs are diamond-free chordal graphs~\cite{bandelt1986}. Linear-time algorithms to compute all three code numbers in block graphs have been presented in~\cite{argiroffo2020linear}. In this paper, we complement these results by determining tight lower and upper bounds for all three code numbers for block graphs. We give bounds using (i) the number of vertices, i.e. the order of a graph, as has been done for several other classes of graphs; and (ii) the number of blocks of a block graph, a quantity equally relevant to block graphs. In doing so, we also prove the following conjecture.

\begin{conjecture} [{\cite[Conjecture 1]{ABLW_ICGT}}]\label{Conj_ID}
The ID-number of a closed-twin-free block graph is bounded above by the number of blocks in the graph.
\end{conjecture}

In addressing LD-codes for twin-free block graphs, we prove (for block graphs) the following conjecture posed by Garijo et al.~\cite{garijo2014difference} and reformulated in a slightly stronger form by Foucaud et al. \cite{foucaud2016location}.

\begin{conjecture} [{\cite[Conjecture 2]{foucaud2016location}}] \label{Conj_LD twin free}
Every twin-free graph $G$ with no isolated vertices satisfies $\gamma^{LD} (G) \leq \frac{|V(G)|}{2}$.
\end{conjecture}

A short version of this article has also been published in the conference proceedings of CALDAM 2023 ~\cite{DBLP:conf/caldam/ChakrabortyFPW23}.

\subsection{Notations}

\begin{figure}[ht!]
     \centering
         \includegraphics[scale=0.48]{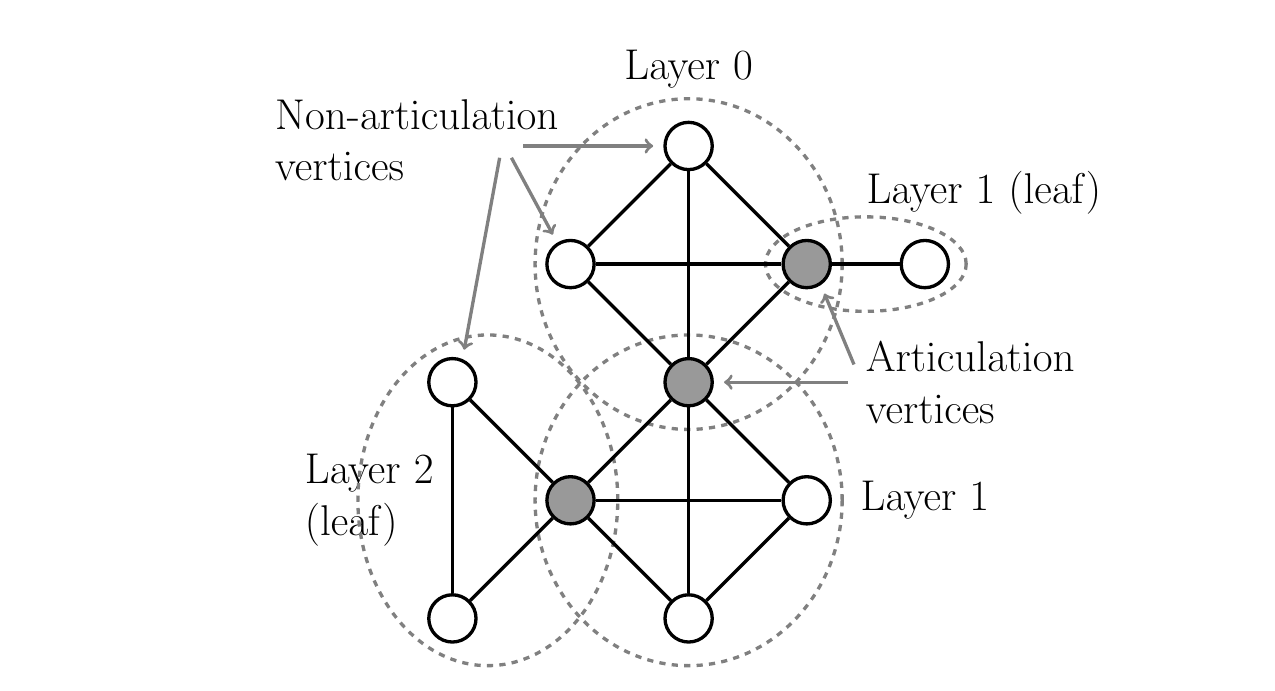}
     \caption{Example of different layer numbers, articulation vertices (grey) and non-articulation vertices (white) of a connected block graph.}
     \label{fig_Layers}\vspace*{-2mm}
\end{figure}

For a block graph $G$, we let $\mathcal{K}(G)$ denote the set of all blocks of $G$, i.e. the set of all maximal cliques of $G$. Noting that any two distinct blocks $K$ and $K'$ of $G$ intersect at at most a single vertex, any vertex $x \in V(G)$ such that $\{ x \}=V(K) \cap V(K')$ is called an \emph{articulation vertex} of both $K$ and $K'$. We define $art(K)$ to be the set of all articulation vertices of a block $K \in \mathcal{K}(G)$. For a connected block graph, we fix a \emph{root block} $K_0 \in \mathcal{K}(G)$ and define a system of assigning numbers to every block of $G$ depending on ``how far" the latter is from $K_0$. So, define a \emph{layer function} $f: \mathcal{K}(G) \to \mathbb{Z}$ \emph{on} $G$ by: $f(K_0)=0$, and for any other $K (\neq K_0) \in \mathcal{K}(G)$ (also called a \emph{non-root} block), define inductively $f(K)=i$ if $V(K) \cap V(K') \neq \emptyset$ for some block $K' ~(\neq K) \in \mathcal{K}(G)$ such that $f(K')=i-1$. For a pair of blocks $K, K' \in \mathcal{K}(G)$ such that $f(K)=f(K')+1$, define $art^-(K)=V(K) \cap V(K')$; and for the root block $K_0$, define $art^- (K_0)=\emptyset$. Note that for a block $K \in \mathcal{K}(G)$ such that $f(K) \geq 1$, we have $|art^-(K)|=1$, and the only vertex in $art^-(K)$ is called the \emph{negative articulation} vertex of the block $K$. In contrast to the negative articulation vertices of $G$, define $art^+ (K) = art(K) \setminus art^- (K)$ to be the set of all \emph{positive articulation} vertices of the block $K$ and $\overline{art} (K) = V(K) \setminus art(K)$ to be the set of all \emph{non-articulation} vertices of $K$. Any block $K$ with $|art(K)|=1$ is called a \emph{leaf block} and all blocks that are not leaf blocks are called \emph{non-leaf blocks}. Note that a block graph always has leaf-blocks. For simplicity, we also denote the set $f^{-1}(\{ i \})$ by $f^{-1}(i)$. Then, for each $i \geq 0$, $f^{-1}(i)$ is called the \emph{$i$-th layer} of $G$ and each block $K \in f^{-1}(i)$ is said to be \emph{in the} $i$-th layer of $G$. See Figure~\ref{fig_Layers} for an illustration of the layers and the related concepts in a connected block graph.

\subsection{Structure of the paper}

Our results on the upper bounds for the code numbers are contained in Section~\ref{sec:UB} of this paper, whereas Section~\ref{sec:LB} is dedicated to the lower bounds for the code numbers. Section~\ref{sec:UB} is further subdivided into three subsections with each of the latter containing the results for a particular code. We conclude in Section~\ref{sec:conclu}.

\section{Upper bounds}\label{sec:UB}

In this section, we establish upper bounds on the ID-, LD- and OLD-numbers for block graphs. Two of these upper bounds are in fact proving Conjectures \ref{Conj_ID} and \ref{Conj_LD twin free}. All our results in this section are for \emph{connected} block graphs. However, applying the results to each connected component of a block graph, the results of Theorem \ref{thm_ub}, \ref{th_LD_ub} and \ref{thm_twin-free} hold equally well for all block graphs.

\subsection{Identifying codes}

The number of blocks is, structurally speaking, a quantity as relevant for block graphs as is the number of vertices for trees. In the next result, we prove Conjecture \ref{Conj_ID} to provide an upper bound for $\gamma^{ID}(G)$ for a block graph $G$ in terms of its number of blocks.

\begin{theorem}\label{thm_ub}
Let $G$ be a connected closed-twin-free block graph and let $\mathcal{K}(G)$ be the set of all blocks of $G$. Then $\gamma^{ID}(G)\leq |\mathcal{K}(G)|$.
\end{theorem}

\begin{proof}
Assume by contradiction that there is a closed-twin-free block graph $G$ of minimum order such that $\gamma^{ID}(G)> |\mathcal{K}(G)|$. We also assume that $G$ has at least four vertices since it can be easily checked that the theorem is true for closed-twin-free block graphs with at most three vertices (which are only $P_1$ and $P_3$). Suppose that $K \in \mathcal{K}(G)$ is a leaf-block of $G$. Due to the closed-twin-free property of $G$, one can assume that $V(K) = \{ x,y \}$ and, without loss of generality, that $x$ and $y$ are the non-articulation and the negative articulation vertices, respectively, of $K$.  Let $G'=G-x$ be the graph obtained by deleting the vertex $x \in V(G)$ (and the edge incident on $x$) from $G$. Then $G'$ is a block graph with $|\mathcal{K}(G')| = |\mathcal{K}(G)| - 1$. We now consider the following two cases.

\begin{case}[$G'$ is closed-twin-free]
By the minimality of the order of $G$, there is an ID-code $C'$ of $G'$ such that $|C'| \leq |\mathcal{K}(G')| = |\mathcal{K}(G)|-1$.
First, assume that $y \notin C'$. Then by the property of domination of $C'$, there exists a vertex $z \in V(G')$ such that $z \in N_{G'}(y) \cap C'$. We claim that $C=C'\cup \{x\}$ is an ID-code of $G$. First of all, that $C$ is a dominating set of $G$ is clear from the fact that $C'$ is a dominating set of $G'$. To prove that $C$ is a closed-separating set of $G$, we see that $x$ is closed-separated in $G$ from all vertices in $V(G') \setminus \{ y \}$ by itself and is closed-separated in $G$ from $y$ by the vertex $z \in C'$. Moreover, all other pairs of distinct vertices are closed-separated by $C'$ are also closed-separated by $C$. Thus, $C$, indeed, is an ID-code of $G$. This implies that $\gamma^{ID}(G) \leq |C| \leq |\mathcal{K}(G)|$, contrary to our assumption.

\smallskip
We therefore assume that $y\in C'$. If again, there exists a vertex $z \in N_{G'}(y) \cap C'$, then by the same reasoning as above, $C = C'\cup\{x\}$ is an ID-code of $G$. Otherwise, we have $N[y]\cap C'=\{y\}$. Now, since $G$ is connected, we have $deg_G(y) > 1$ and therefore, there exists a vertex $w \in N_G(y) \setminus \{ x \}$. Then $C = C' \cup \{ w \}$ is an ID-code of $G$. This is because, first of all, $C$ still closed-separates every pair of distinct vertices in $V(G')$. The vertex $x$ is closed-separated from $y$ by the vertex $w \in C$; and from $w$ by $w$ itself. Moreover, for any vertex $v$ in $V(G') \setminus \{ y,w \}$, $v$ is closed-separated from $y$ in $G'$ by some vertex $u_v \in N_{G'}[v] \cap C'$. Then $x$ is closed-separated from all such $v$ in $V(G') \setminus \{ y,w \}$ by the vertices $u_v \in C'$. Moreover, $C$ is clearly also a dominating set of $G$. Hence, this leads to the same contradiction as before.
\end{case}

\begin{case}[$G'$ has closed-twins]
Assume that vertices $u$, $v \in V(G')$ are a pair of closed-twins of $G'$. Since $u$ and $v$ were not closed-twins in $G$, it means that $x$ is adjacent to, say, $u$, without loss of generality. This implies that $u=y$. Note that $v$ is then unique with respect to being a closed-twin with $y$ in $G'$. This is because, if $u$ and some vertex $v' (\neq v) \in V(G')$ were also closed-twins in $G'$, then it would mean that $v$ and $v'$ were closed-twins in $G$, contrary to our assumption. Now, let $G''=G'- v$. We claim the following.

\begin{caseclaim}
$G''$ is closed-twin-free.
\end{caseclaim}

\begin{proofofcaseclaim}
Toward a contradiction, if vertices $z,w \in V(G'')$ were a pair of closed-twins in $G''$, it would then mean that the vertex $z \in N_{G'}(v)$, without loss of generality. This would, in turn, imply that $z \in N_{G'}(y)$ (since the vertices $y$ and $v$ are closed-twins in $G'$). Or, in other words, $y \in N_{G''}(z)$. Now, since $z$ and $w$ are closed-twins in $G''$, we have $y \in N_{G''}(w)$, i.e. $w \in N_{G'}(y)$. Again, by virtue of $y$ and $v$ being closed-twins in $G'$, we have $w \in N_{G'}(v)$. This implies that $z$ and $w$ are closed-twins in $G$ which is a contradiction to our assumption.
\end{proofofcaseclaim}

We also note here that the vertices $y$ and $v$ must be from the same block, as the two are adjacent on account of being closed-twins in $G'$. Thus, $G''$ is a connected closed-twin-free block graph. Therefore, by the minimality of the order of $G$, there is an ID-code $C''$ of $G''$ such that $|C''| \leq |\mathcal{K}(G'')| < |\mathcal{K}(G)|$. If $y\notin C''$, then we claim that $C = C''\cup \{x\}$ is an identifying code of $G$. This is true because, firstly, $C$ is a dominating set of $G$ (note that, by the property of domination of $C''$ in $G''$, there exists a vertex $z \in N_{G''}(y) \cap C''$; and since $y$ and $v$ are closed-twins in $G'$, we have $z \in N_{G}(v) \cap C$). Moreover, $x$ is closed-separated in $G$ from every other vertex in $V(G) \setminus \{ y \}$ by $x$ itself; and the vertices $x$ and $y$ are closed-separated in $G$ by some vertex in $N_{G''}(y) \cap C''$ that dominates the vertex $y$. The vertex $y$ is closed-separated from all the vertices in $V(G) \setminus \{y,x\}$ by $x$ and; since $y$ and $v$ have the same closed neighborhood in $G'$, $v$ is closed-separated in $G$ from all vertices in $V(G'') \setminus \{ y \}$ because $y$ is by $C''$. Finally, any two distinct vertices closed-separated by $C''$ still remain so by $C$. Thus, $C$, indeed, is an ID-code of $G$. This implies that $\gamma^{ID}(G) \leq |C| \leq |\mathcal{K}(G)|$; again a contradiction.

\smallskip
Let us, therefore, assume that $y\in C''$. This time, we claim that $C=(C'' \setminus \{ y \}) \cup\{x,v\}$ is an ID-code of $G$. That $C$ is a dominating set of $G$ is clear. So, as for the closed-separating property of $C$ is concerned, as before, $x$ is closed-separated in $G$ from every other vertex in $V(G) \setminus \{ y \}$ by $x$ itself; vertices $x$ and $y$ are closed-separated in $G$ by $v$; the vertices $y$ and $v$ are closed-separated in $G$ by $x$ and the vertices $v$ and $x$ are closed-separated in $G$ by $v$. Since $y$ and $v$ have the same closed neighbourhood in $G'$ and since $y$ is closed-separated in $G''$ by $C''$ from every other vertex in $V(G'')$, both $v$ and $y$ are also each closed-separated in $G$ from every vertex in $V(G'') \setminus \{ v,y \}$.

Finally, any two distinct vertices of $G''$ are closed-separated by $C''$ still remain so by $C$. This proves that $C$ is an ID-code of $G$ and hence, again, we are led to the contradiction that $\gamma^{ID}(G) \leq |C| \leq |\mathcal{K}(G)|$.
\end{case}

This proves the theorem.
\end{proof}

Note that, besides for stars \cite{DBLP:journals/dm/GravierM07}, the upper bound in Theorem \ref{thm_ub} is attained by the ID-numbers of thin headless spiders~\cite{argiroffo2014study} which, therefore, serve as examples of cases where the bound in Theorem \ref{thm_ub} is tight. Noting that blocks are simply the maximal complete subgraphs of a block graph, the statement of Theorem \ref{thm_ub} does not hold in terms of the number of maximal complete subgraphs even for chordal graphs, let alone for graphs in general. A counterexample to the bound in Theorem \ref{thm_ub} (replacing the number of blocks by the number of maximal cliques) for chordal graphs is the graph $P_{2k}^{k-1}$ (the graph obtained from a path on $2k$ vertices with edges introduced between all pairs of vertices $u,v \in V(P_{2k})$ such that $d_{P_{2k}}(u,v) \leq k-1$) which is closed-twin-free, has only two maximal complete subgraphs, but needs $2k-1$ vertices in any identifying code \cite{foucaud2011extremal}.

\subsection{Locating-dominating codes}

In this subsection, we prove two results on upper bounds for the LD-numbers of block graphs. The first result is a more general one in which the upper bound is in terms of the number of blocks and other quantities arising out of the structural properties of a block graph. On the other hand, the second result is proving Conjecture \ref{Conj_LD twin free} for block graphs. We begin with the more general result.

\begin{theorem} \label{th_LD_ub}
Let $G$ be a connected block graph and $\mathcal{K}(G)$ be the set of blocks in $G$. Then, we have
$$\gamma^{LD}(G) \leq |\mathcal{K}(G)| + \sum_{\substack{K \in \mathcal{K}(G), \\|\overline{art} (K)| \geq 2}} (|\overline{art} (K)|-2).$$
\end{theorem}

\begin{proof}
We define a set $C \subset V(G)$ by the following rules.

\begin{enumerate}[{Rule} 1:]
	\item For every block $K \in \mathcal{K}(G)$ which does not contain any closed-twins, i.e. with at most one non-articulation vertex, pick any one vertex from $V(K) \setminus art^-(K)$ in $C$.
	
	\item For every block $K \in \mathcal{K}(G)$ which contains closed-twins, i.e. with at least two non-articulation vertices, pick any $|\overline{art} (K)|-1$ vertices from $\overline{art} (K)$ in $C$.
\end{enumerate}

Note that the vertices added in $C$ by the above rules are all distinct. Therefore, the following is the size of $C$.
\begin{flalign*}
|C| &= \Big|\{ K \in \mathcal{K}(G): |\overline{art} (K)| \leq 1 \} \Big| + \sum_{\substack{K \in \mathcal{K}(G),\\ |\overline{art} (K)| \geq 2}} (|\overline{art} (K)|-1)\\
&= |\mathcal{K}(G)| + \sum_{\substack{K \in \mathcal{K}(G), \\|\overline{art} (K)| \geq 2}} (|\overline{art} (K)|-2).
\end{flalign*}
The result, therefore, follows from proving that $C$ is an LD-code of $G$.

\medskip
First of all, we notice that, by the construction of $C$, for every block $K \in \mathcal{K}(G)$, there exists a vertex $v_K \in V(K) \cap C$. Therefore, $C$ is a dominating set of $G$. We now show that $C$ is also a locating set of $G$. So assume that $u,v \in V(G) \setminus C$ are distinct vertices of $G$. Then, by the construction of $C$, we must have $u \in V(K)$ and $v \in V(K')$ for a distinct pair of blocks $K, K' \in \mathcal{K}(G)$.

\begin{claim}
There exist vertices $v_K \in V(K) \cap C$ and $v_{K'} \in V(K') \cap C$ such that $v_K \neq v_{K'}$.
\end{claim}

\begin{proofofclaim}
On the contrary, if $V(K) \cap C  = V(K') \cap C = \{v_K\}$, this would imply that $v_K$ is the negative articulation vertex of either $K$ or $K'$. Without loss of generality, let us assume that $art^-(K) = \{v_K\}$. If $K$ does not contain any closed twins, then by Rule 1, there exists a vertex in $K$ other than $v_K$ which belongs to $C$, a contradiction to $V(K) \cap C = \{v_K\}$. Therefore, $K$ contains closed twins. In other words, $|\overline{art}(K)| \geq 2$. Now, by Rule 2, at least one vertex $w$, say, of $\overline{art}(K)$ belongs to $C$. Since $\overline{art}(K) \cap art^-(K) = \emptyset$, we have $w \neq v_K$ and thus, we run into the same contradiction. This proves the claim.
\end{proofofclaim}

This implies that at least one of $v_K$ and $v_{K'}$ must open-separate $u$ and $v$ in $G$ and, hence, $C$ is a locating set of $G$.
\end{proof}

\begin{figure}[h!]
\vspace*{-2mm}
     \centering
         \includegraphics[scale=0.4]{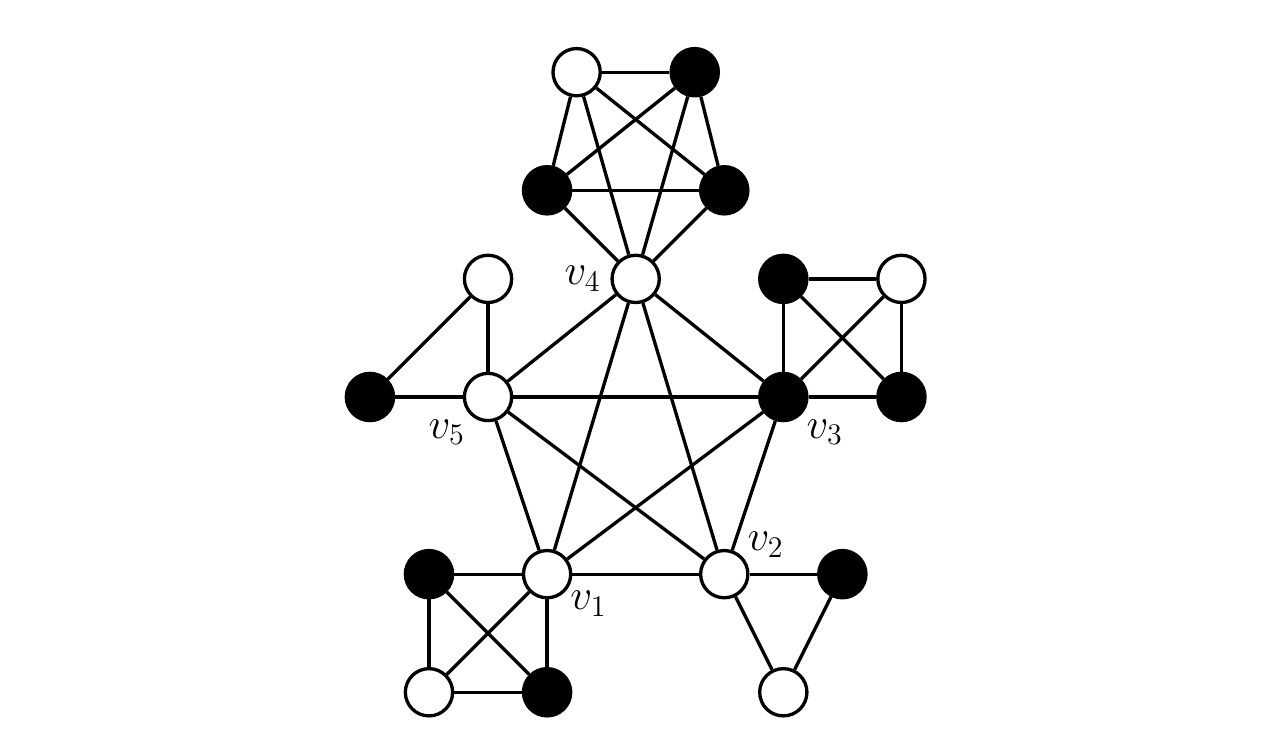}\vspace*{-1mm}
     \caption{Graph $S_5 (4,3,4,5,3)$ whose LD-number attains the upper bound in Theorem \ref{th_LD_ub}. The black vertices represent those included in the LD-code $C$ of $G$ as described in the proof of Theorem \ref{th_LD_ub}.}
     \label{fig_LD_ub}\vspace*{-2mm}
\end{figure}

There is an infinite number of arbitrarily large connected block graphs whose LD-numbers attain the upper bound in Theorem \ref{th_LD_ub}. One such subclass of block graphs is the following. For positive integers $t \geq 2$ and $m_1, m_2, \ldots, m_t$ with $m_i \geq 3$ for each $i$, we define a class of graphs $S_t(m_1, m_2, \ldots , m_t)$ by the following rule. Let $X$ be a copy of the complete graph on $t$ vertices and name its vertices $v_1, v_2, \ldots , v_t$. Also, for all $1 \leq i \leq t$, let $Y_i$ be a copy of the complete graph on $m_i$ vertices. Let $S_t(m_1, m_2, \ldots , m_t)$ be the block graph obtained by identifying a vertex of $Y_i$ with $v_i$ of $X$ for every $1 \leq i \leq t$. For brevity, we continue to call the identified vertices resulting in $S_t(m_1, m_2, \ldots , m_t)$ by the same names of $v_1, v_2, \ldots , v_t$ as before. See Figure \ref{fig_LD_ub} for an example of the graph $S_t(m_1, m_2, \ldots , m_t)$ constructed with $t=5$, and $m_1=m_3=4$, $m_2=m_5=3$ and $m_4=5$. We then show the following.

\begin{proposition} \label{prop}
For $t \geq 2$, and $m_1, m_2, \ldots , m_t$ such that $m_i \geq 3$ for all $1 \leq i \leq t$, we have 
$$\gamma^{LD}(S_t (m_1, m_2, \ldots , m_t)) = |\mathcal{K}(S_t (m_1, m_2, \ldots , m_t))| + \sum_{\substack{K \in \mathcal{K}(S_t (m_1, m_2, \ldots , m_t)), \\|\overline{art} (K)| \geq 2}} (|\overline{art} (K)| - 2).$$
\end{proposition}

\begin{proof}
Let $G = S_t (m_1, m_2, \ldots , m_t)$. We note here that the number of blocks in $G$ is $t+1$; and the only blocks $K \in \mathcal{K}(G)$ with $|\overline{art} (K)| \geq 2$ are $Y_1, Y_2, \ldots , Y_t$ (as per the notations used in the preceeding discussion). More precisely, for each $1 \leq i \leq t$, we have $|\overline{art} (Y_i)| = m_i-1$. So, the upper bound for the LD-number of $G$ by Theorem \ref{th_LD_ub} is $t+1+ \sum_{i=1}^t (m_i-3) = 1-2t + \sum_{i=1}^t m_i$.

Now, assume that $C$ is a minimum LD-code of $G$. Let us first assume that $V(X) \cap C = \emptyset$. Then, since any two vertices in $V(Y_i)$ have the same neighborhood in $C$, it implies that we must have $|V(Y_i) \cap C| = m_i-1$. This further implies that $|C| = \sum_{i=1}^t (m_i-1) = -t + \sum_{i=1}^t m_i > 1-2t + \sum_{i=1}^t m_i$ (since $t \geq 2$), the upper bound by Theorem \ref{th_LD_ub} resulting in a contradiction. Therefore, we must have $V(X) \cap C \neq \emptyset$. Thus, let $v_i \in C$ for some $1 \leq i \leq t$. Note that $v_i \in V(Y_i)$. Since any two vertices of $V(Y_i) \setminus \{v_i\}$ are twins in $G$, we have $|V(Y_i) \cap C| \geq m_i-1$. Moreover, for $1 \leq j \leq t$ such that $j \neq i$, again, since any two vertices of $V(Y_j) \setminus \{v_j\}$ are twins in $G$, we now have $|V(Y_j) \cap C| \geq m_j - 2$. Hence, we have $|C| \geq |V(Y_i) \cap C| + \sum_{j=1, j \neq i}^t |V(Y_j) \cap C| \geq 1-2t + \sum_{i=1}^t m_i$.
\end{proof}

Imposing additional structural constraints on a block graph, one could still limit the number of vertices one needs to choose from each of its blocks in order to form an LD-code of the graph. Our next result shows exactly that.

\begin{theorem}\label{thm_twin-free}
Let $G$ be a connected twin-free block graph. Then $\gamma^{LD}(G) \leq \frac{|V(G)|}{2}$.
\end{theorem}

\begin{proof}
To prove the theorem, we partition the vertex set of $G$ into two special subsets $C^*$ and $D^*$.

\medskip
Assume that $K_0 \in \mathcal{K}(G)$ is a leaf block of $G$. Then, $|V(K_0)|=2$, as $G$ is twin-free. Assign $K_0$ to be the root block of $G$, i.e. define a layer function $f : \mathcal{K}(G) \to \mathbb{Z}$  on $G$ such that $f(K_0) = 0$. We then construct the sets $C^*$ and $D^*$ by the following rules applied inductively on $i \in f(\mathcal{K}(G))$. See Figure~\ref{constr_fig} for a demonstration of this construction.

\begin{figure}[!h]
\vspace*{3mm}
     \centering
     \begin{subfigure}[!h]{0.19\textwidth}
         \centering
    \hspace*{-6mm}     \includegraphics[width=\textwidth]{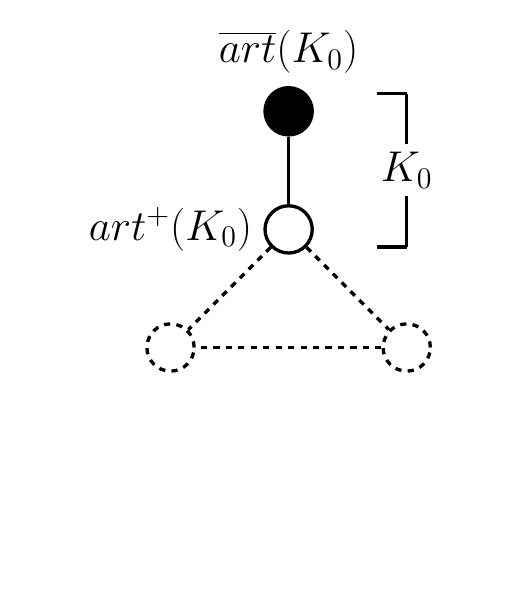}
         \caption{Rule 1 }
         \label{Rule 1}
     \end{subfigure}
     \begin{subfigure}[!h]{0.19\textwidth}
         \centering
         \includegraphics[width=\textwidth]{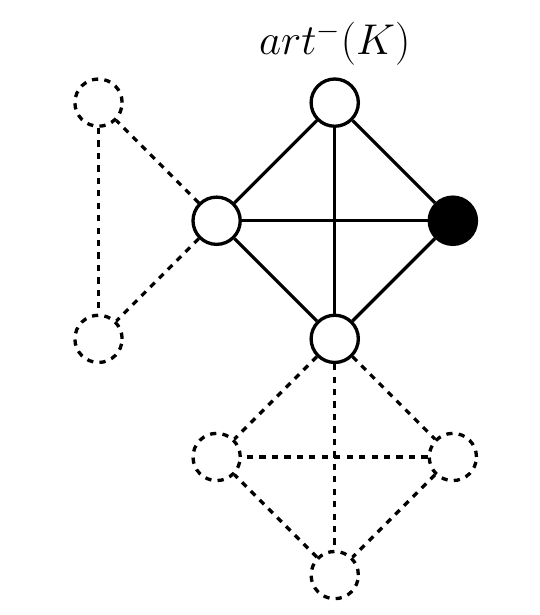}
         \caption{Rule 2}
         \label{Rule 2}
     \end{subfigure}
     \begin{subfigure}[!h]{0.19\textwidth}
         \centering
         \includegraphics[width=\textwidth]{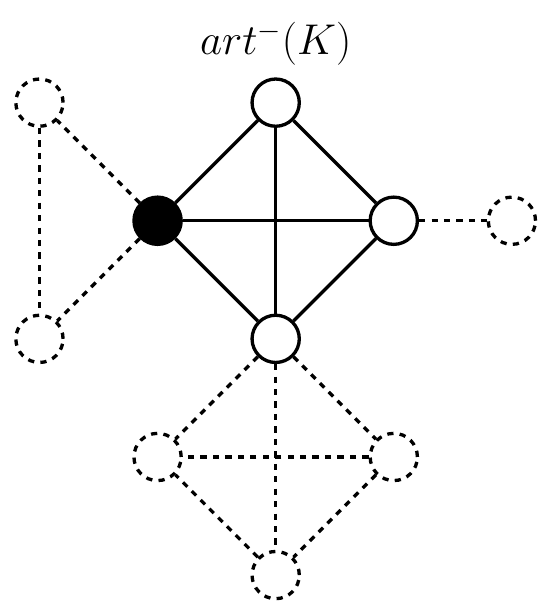}
         \caption{Rule 3}
         \label{Rule 3}
     \end{subfigure}
     \begin{subfigure}[!h]{0.19\textwidth}
         \centering
         \includegraphics[width=\textwidth]{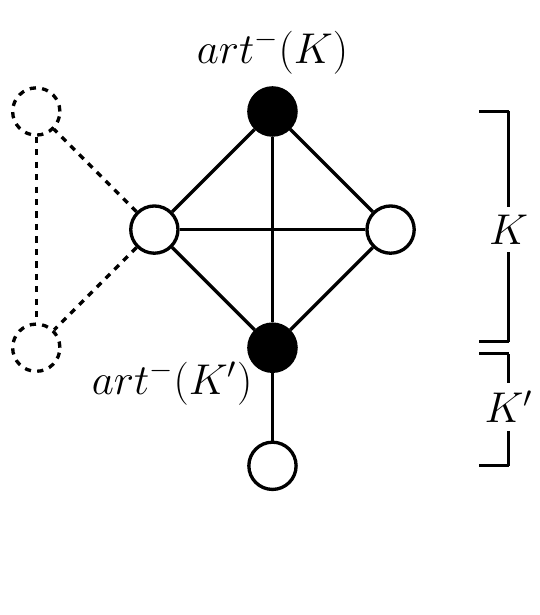}
         \caption{Rule 4}
         \label{Rule 4}
     \end{subfigure}
     \begin{subfigure}[!h]{0.19\textwidth}
         \centering
         \includegraphics[width=\textwidth]{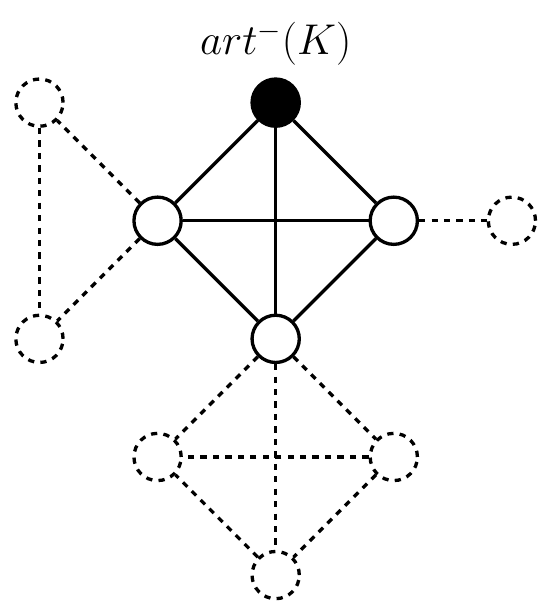}
         \caption{Rule 5}
         \label{Rule 5}
     \end{subfigure}\vspace*{-2mm}
     \caption{Example of each rule in the proof of Theorem \ref{thm_twin-free} for the construction of the sets $C^*$ or $D^*$. In each example, the black vertices represent those picked in $C^*$ and the white vertices represent those picked in $D^*$. The blocks with solid edges represent those blocks (in $i$-th layer, say) from which vertices are chosen either in $C^*$ or $D^*$. The blocks with dashed edges represent those blocks in the next layer (the $(i+1)$-th) which, inductively, are yet to be analysed for their choices of vertices in $C^*$ and $D^*$; but whose presence in the figure is necessary to determine the positive, negative and the non-articulation vertices of the block in $i$-th layer.}
     \label{constr_fig}\vspace*{-4mm}
\end{figure}

\begin{enumerate}[{Rule} 1:]
	\item Pick the (positive) articulation vertex of the root block $K_0$ in $D^*$ (i.e. let $art^+ (K_0) \subset D^*$) and pick the (other) non-articulation vertex of $K_0$ in $C^*$ (i.e. let $\overline{art} (K_0) \subset C^*$). See Figure \ref{Rule 1} for an example.
	
	\item For every non-root block $K \in \mathcal{K}(G)$ with at least one non-articulation vertex (i.e. $\overline{art} (K) \neq \emptyset$) and whose negative articulation vertex is in $D^*$ (i.e. $art^- (K) \subset D^*$), pick all non-articulation vertices of $K$ in $C^*$ (i.e. let $\overline{art} (K) \subset C^*$); and all positive articulation vertices of $K$ in $D^*$ (i.e. let $art^+ (K) \subset D^*$). See Figure \ref{Rule 2} for an example.
	
	\item For every non-root block $K \in \mathcal{K}(G)$ with no non-articulation vertices (i.e. $\overline{art} (K) = \emptyset$) and whose negative articulation vertex is in $D^*$ (i.e. $art^- (K) \subset D^*$), pick \emph{one} positive articulation vertex, say, $w$ of $K$ in $C^*$ and the rest of the positive articulation vertices in $D^*$ (i.e. let $art^+ (K) \setminus \{ w \}\subset D^*$). See Figure \ref{Rule 3} for one such case.
	
	\item For every non-root block $K \in \mathcal{K}(G)$ with at least one non-articulation vertex (i.e. $\overline{art} (K) \neq \emptyset$) and whose negative articulation vertex is in $C^*$ (i.e. $art^- (K) \subset C^*$), pick \emph{one} positive articulation vertex (if available), say, $w$ of $K$ in $C^*$; and pick all other vertices in $V(K)$, except the vertex $w$ and the negative articulation vertex of $K$, in $D^*$ (i.e. let $V(K) \setminus (art^- (K) \cup \{ w \}) \subset D^*$). See Figure \ref{Rule 4} for both examples of when articulation vertices are available (block $K$) and when they are not, i.e. in the case of leaf blocks (block $K'$).
	
	\item For every non-root block $K \in \mathcal{K}(G)$ with no non-articulation vertices (i.e. $\overline{art} (K) = \emptyset$) and whose negative articulation vertex is in $C^*$ (i.e. $art^- (K) \subset C^*$), pick all positive articulation vertices of $K$ in $D^*$ (i.e. $art^+ (K) \subset D^*$). See Figure \ref{Rule 5} for an illustration.
\end{enumerate}

From the construction, $C^*$ and $D^*$ are complements of each other in $V(G)$. We claim that both $C^*$ and $D^*$ are LD-codes of $G$. We first show that both are dominating sets of $G$.

\begin{claim} \label{clm_twin-free_dom}
Both $C^*$ and $D^*$ are dominating sets of $G$.
\end{claim}

\begin{proofofclaim}
To prove that both $C^*$ and $D^*$ are dominating sets of $G$, it is enough to show that, for every block $K \in \mathcal{K}(G)$, both $V(K) \cap C^* \neq \emptyset$ and $V(K) \cap D^* \neq \emptyset$. By Rule 1, the claim is true for the root block $K_0$. So, assume $K \in \mathcal{K}(G)$ to be a non-root block. First, suppose that the negative articulation vertex of $K$ belongs to $D^*$. Then, by Rules 2 and 3, we have $V(K) \cap C^* \neq \emptyset$. Next, suppose that the negative articulation vertex of $K$ belongs to $C^*$. Then, by Rules 4 and 5, we have $V(K) \cap D^* \neq \emptyset$.
\end{proofofclaim}

We now show that both $C^*$ and $D^*$ are also locating sets of $G$. We start with $C^*$.

\begin{claim} \label{clm_twin-free_C*}
$C^*$ is a locating set of $G$.
\end{claim}

\begin{proofofclaim}
Assume that $u, v \in D^*$ are distinct vertices of $G$. Since $G$ is twin-free, there exist distinct blocks $K, K' \in \mathcal{K}(G)$ such that $u \in V(K)$ and $v \in V(K')$. By the proof of Claim \ref{clm_twin-free_dom}, there exist vertices $v_K \in V(K) \cap C^*$ and $v_{K'} \in V(K') \cap C^*$. If $v_K \neq v_{K'}$, then either one of $v_K$ and $v_{K'}$ must locate $u$ and $v$ in $G$. So, let us assume that no such pairs of distinct vertices $v_K \in V(K) \cap C^*$ and $v_{K'} \in V(K') \cap C^*$ exist, i.e. $V(K) \cap V(K') \subset C^*$ and that $V(K) \triangle V(K') \subset D^*$.

\medskip
We now claim that either $u$ is an articulation vertex of $K$ or $v$ is an articulation vertex of $K'$ (or both). So, toward contradiction, assume that both $u$ and $v$ are non-articulation vertices of $K$ and $K'$, respectively. Then the following two cases arise.

\begin{case}[$K$ and $K'$ belong to different layers]
Without loss of generality, assume that $f(K')=f(K)+1$. Then, $K \neq K_0$, or else, by Rule 1, $u$, being a non-articulation vertex of $K$, must belong to $C^*$, contrary to our assumption. Therefore, $K$ is a non-leaf block. Now, $V(K) \triangle V(K') \subset D^*$ implies that the negative articulation vertex of $K$ belongs to $D^*$. Since $u$ is a non-articulation vertex of $K$, by Rule 2, $u \in C^*$ which is a contradiction to our assumption.
\end{case}

\begin{case}[$K$ and $K'$ belong to the same layer]
In this case, $K$ and $K'$ cannot both be leaf blocks, or else, $G$ would have twins. So, without loss of generality, suppose that $K$ is a non-leaf block. Now, the negative articulation vertex of $K$ belongs to $C^*$. Since $K$ is a non-leaf block, there exists a positive articulation vertex of $K$ and, hence, by Rule 4, $art^+ (K) \cap C^* \neq \emptyset$ which contradicts the fact that $V(K) \triangle V(K') \subset D^*$.
\end{case}

This proves our claim that either $u$ is an articulation vertex of $K$ or $v$ is an articulation vertex of $K'$ (or both). So if, without loss of generality, we assume that $u$ is an articulation vertex of $K$, then $\{ u \} = V(K) \cap V(K'')$ for some block $K'' ~(\neq K) \in \mathcal{K}(G)$. Moreover, $K'' \neq K'$, or else, $V(K) \cap V(K') = \{u\} \subset D^*$ which contradicts our assumption that $V(K) \cap V(K') \subset C^*$. Hence, some vertex in $V(K'') \cap C^*$ (which exists due to the proof of Claim \ref{clm_twin-free_dom}) must open-separate $u$ and $v$ in $G$. This proves our current claim.
\end{proofofclaim}

We now prove the same for $D^*$.

\begin{claim} \label{clm_twin-free_D*}
$D^*$ is a locating set of $G$.
\end{claim}

\begin{proofofclaim}
Assume that $u, v \in C^*$ are distinct vertices of $G$. Since $G$ is twin-free, there exist distinct blocks $K, K' \in \mathcal{K}(G)$ such that $u \in V(K)$ and $v \in V(K')$. By the proof of Claim \ref{clm_twin-free_dom}, there exist vertices $v_K \in V(K) \cap D^*$ and $v_{K'} \in V(K') \cap D^*$. If $v_K \neq v_{K'}$, then either one of $v_K$ and $v_{K'}$ must open-separate $u$ and $v$ in $G$. So, let us assume that no such pairs of distinct vertices $v_K \in V(K) \cap D^*$ and $v_{K'} \in V(K') \cap D^*$ exist, i.e. $V(K) \cap V(K') \subset D^*$ and that $V(K) \triangle V(K') \subset C^*$.

\medskip
We now claim that either $u$ is an articulation vertex of $K$ or $v$ is an articulation vertex of $K'$ (or both). So, toward contradiction, assume that both $u$ and $v$ are non-articulation vertices of $K$ and $K'$, respectively. Then the following two cases arise.

\setcounter{case}{0}
\begin{case}[$K$ and $K'$ belong to different layers]
Without loss of generality, assume that $f(K')=f(K)+1$. If $|V(K')| \geq 3$, since $G$ is twin-free and since $v$ is a non-articulation vertex of $K'$, then $K'$ contains exactly one non-articulation vertex and thus, $art^+ (K') \cap D^* \neq \emptyset$ by Rule 2. This, however, is a contradiction to the fact that $V(K) \triangle V(K') \subset C^*$. So, assume that $|V(K')|=2$, in which case, $K'$ is a leaf block (since, again, $v$ is a non-articulation vertex of $K'$). This implies that $K$ is a non-leaf block, or else, $G$ would have twins. So, in particular, $K \neq K_0$, the root block of $G$. Moreover, $V(K) \triangle V(K') \subset C^*$ implies that the negative articulation vertex of $K$ belongs to $C^*$. Therefore, since $u$ is a non-articulation vertex of $K$, by Rule 4, $u \in D^*$ which is a contradiction to our assumption.
\end{case}

\begin{case}[$K$ and $K'$ belong to the same layer]
In this case, $K$ and $K'$ cannot both be leaf blocks, or else, $G$ would have twins. So, without loss of generality, assume $K$ to be a non-leaf block. Therefore, $|V(K)| \geq 3$, or else, $u$ would be an articulation vertex of $K$, contrary to our assumption. The negative articulation vertex of $K$ belongs to $D^*$. Therefore, by Rule 2, $art^+ (K) \cap D^* \neq \emptyset$ which contradicts $V(K) \triangle V(K') \subset C^*$.
\end{case}

This, therefore, proves our claim that either $u$ is an articulation vertex of $K$ or $v$ is an articulation vertex of $K'$ (or both). If, without loss of generality, $u$ is an articulation vertex of $K$, then $\{ u \} = V(K) \cap V(K'')$ for some block $K'' ~(\neq K) \in \mathcal{K}(G)$. Moreover, $K'' \neq K'$, or else, $V(K) \cap V(K') = \{u\} \subset C^*$ which contradicts our assumption that $V(K) \cap V(K') \subset D^*$. Hence, some vertex in $V(K'') \cap D^*$ (which exists due to the proof of Claim \ref{clm_twin-free_dom}) must open-separate $u$ and $v$ in $G$. This, again, proves our current claim.
\end{proofofclaim}

Combining Claims \ref{clm_twin-free_dom}, \ref{clm_twin-free_C*} and \ref{clm_twin-free_D*}, we find that $C^*$ and $D^*$ are both LD-codes of the twin-free block graph $G$ with no isolated vertices. Moreover, since $C^*$ and $D^*$ are complements of each other in $V(G)$, at least one of them must have cardinality of at most half the order of $G$. This proves the theorem.
\end{proof}

Theorem \ref{thm_twin-free} therefore proves Conjecture \ref{Conj_LD twin free} for block graphs.

\begin{corollary}
Let $G$ be a twin-free block graph without isolated vertices. Then $\gamma^{LD}(G) \leq \frac{|V(G)|}{2}$.
\end{corollary}

\begin{figure}[b!]
     \centering
         \includegraphics[scale=0.4]{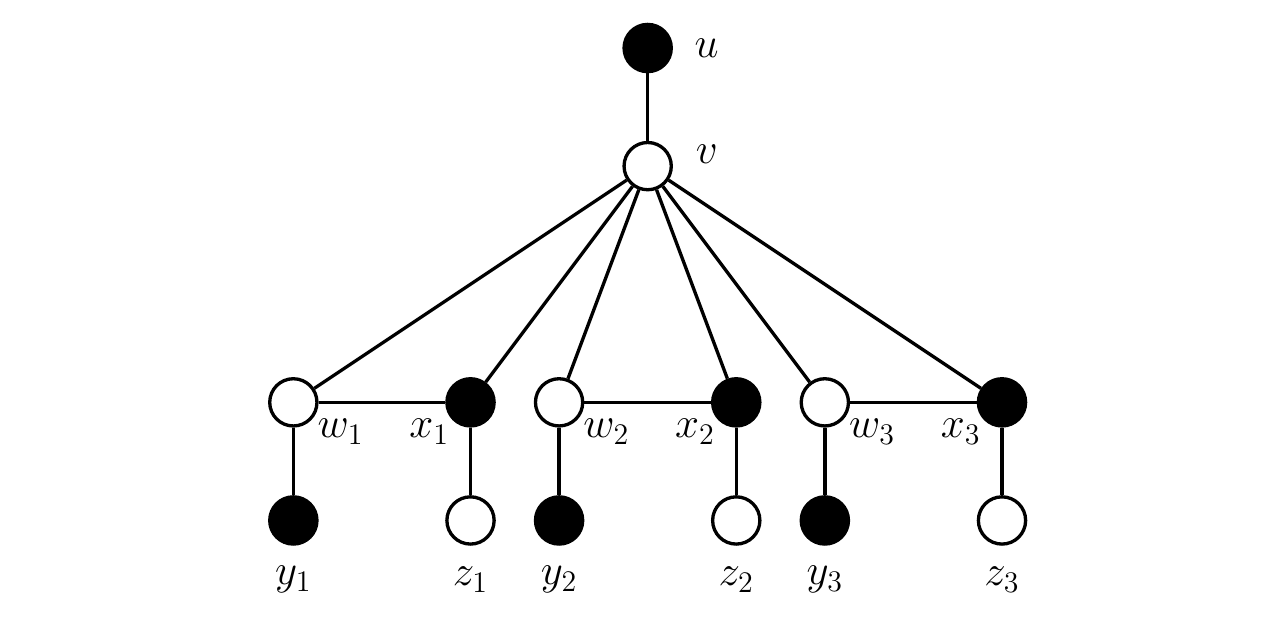}\vspace*{-1mm}
     \caption{Graph $H_3$ whose LD-number attains the upper bound in Theorem \ref{thm_twin-free}. The black vertices represent those included in the LD-code $C^*$ of $G$ described in the proof of Theorem \ref{thm_twin-free}.}
     \label{fig_LD-twin-free}\vspace*{-2mm}
\end{figure}

The trees attaining the bound of Theorem \ref{thm_twin-free} were characterized in \cite{foucaud2016location}. There are also arbitrarily large twin-free block graphs that are not trees and whose LD-numbers attain the bound given in Theorem \ref{thm_twin-free}. To demonstrate this attainment, we look at the following subclass of block graphs which we denote by $H_t$: For a fixed integer $t \geq 1$, let $T_1, T_2, \ldots , T_t$ be $t$ copies of $K_3$, the complete graph on three vertices. Suppose that $V(T_i)= \{ v_i, w_i, x_i \}$ for each $1 \leq i \leq t$. Also, let $R, R_1, R_2, \ldots , R_t, R'_1, R'_2, \ldots , R'_t$ be $2t+1$ copies of $P_2$, the path on two vertices. Also, let $V(R)= \{ u,v \}$ and for all $1 \leq i, i' \leq t$, let $V(R_i) = \{ y'_i, y_i \}$ and $V(R'_{i'}) = \{ z'_{i'}, z_{i'} \}$. We then identify the vertices $v, v_1, v_2, \ldots , v_t$ to a single vertex which we continue to call $v$; and, for each $1 \leq i \leq t$, we identify the vertices $w_i$ and $y'_i$ to a single vertex and the vertices $x_{i}$ and $z'_{i}$ to a single vertex. In the latter two cases, we continue to call the identified vertices $w_i$ and $x_i$, respectively. The new resulting graph is what we call $H_t$. See Figure~\ref{fig_LD-twin-free} for an example of $H_t$ with $t=3$. With that, we now prove the following.

\begin{proposition}
For each integer $t \geq 1$, $\gamma^{LD} (H_t) = \frac{|V(H_t)|}{2}$.
\end{proposition}

\begin{proof}
Notice that the graph $H_t$ is twin-free. Since $|V(H_t)| = 4t+2$, we therefore have from Theorem \ref{thm_twin-free} that $\gamma^{LD} (H_t) \leq 2t+1$.

\medskip
We now prove that $\gamma^{LD} (H_t) \geq 2t+1$. Since each of the $2t+1$ edges $uv$, $w_iy_i$, $x_iz_i$ (for $1\leq i \leq t$) of $H_t$ contains a vertex of degree 1, therefore any LD-code of $H_t$, by its property of domination, must contain at least one endpoint of each of these edges. Since the above edges are all pairwise disjoint, any LD-code of $H_t$ must contain at least $2t+1$ vertices of $H_t$.
\end{proof}

\subsection{Open locating-dominating codes}

\begin{figure}[b!]
\vspace*{-6mm}
     \centering
         \includegraphics[scale=0.4]{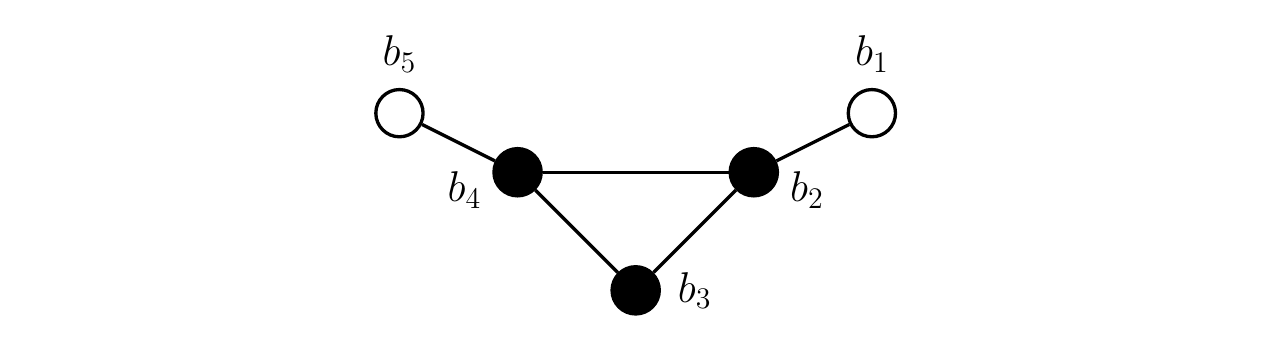}\vspace{-1mm}
     \caption{The Bull graph $B_5$. The set of black vertices constitute an OLD-code of $B_5$.}
     \label{fig_Bull_OLD}

 \vspace*{4mm}
     \centering
     \begin{subfigure}[!h]{0.49\textwidth}
         \centering
         \includegraphics[scale=0.4]{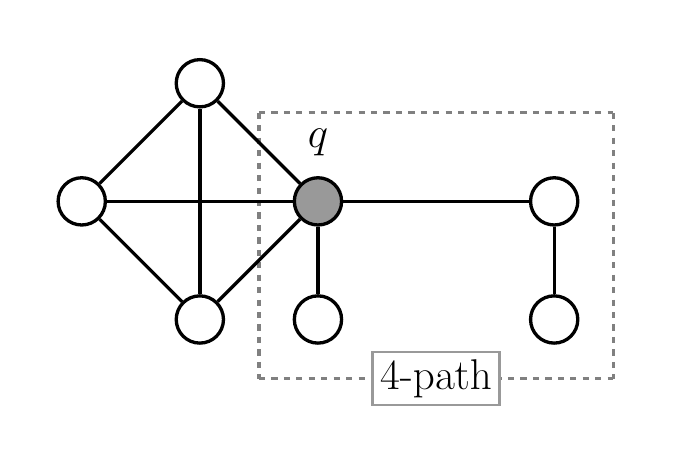}\vspace*{-1mm}
         \caption{$K_4 \rhd_q P_4$}
         \label{fig_G adjunct P4}
     \end{subfigure}
     \begin{subfigure}[!h]{0.49\textwidth}
         \centering
         \includegraphics[scale=0.4]{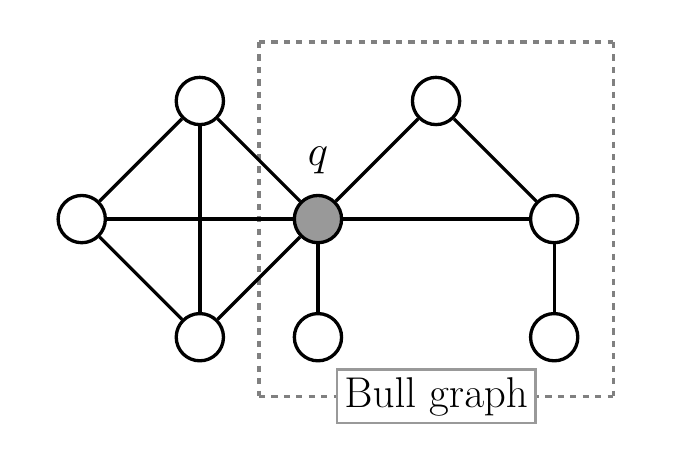}\vspace*{-1mm}
         \caption{$K_4 \rhd_q B_5$}
         \label{fig_G adjunct B5}
     \end{subfigure}\vspace*{-2mm}
     \caption{Examples of $G' \rhd_q X$, where $G' \cong K_4$ and $X \in \{ P_4, B_5 \}$. The vertex $q$ (in grey) is obtained by identifying a vertex of $G'$ and an articulation vertex of $X$.}
     \label{fig_G adjunct X}
\end{figure}

We now focus our attention on upper bounds for OLD-numbers for block graphs. Before we get to our results, we define the following two special graphs.
\begin{enumerate}
	\item The \emph{$4$-path} (or $P_4$ in symbol) is a graph defined by its vertex set $V(P_4)=\{ p_1, p_2, p_3, p_4 \}$ and its edge set $E(P_4)= \{ p_1 p_2, p_2 p_3, p_3 p_4 \}$.
	
	\item The \emph{bull graph} (or $B_5$ in symbol) is a  graph defined by its vertex set $V(B_5)=\{ b_1, b_2, b_3, b_4, b_5 \}$ and its edge set $E(B_5)= \{ b_1 b_2, b_2 b_3, b_3 b_4, b_4 b_5, b_2 b_4 \}$. See Figure \ref{fig_Bull_OLD} for a depiction of a bull graph.
\end{enumerate}

We note here that both $P_4$ and $B_5$ are block graphs with articulation vertices $p_2$ and $p_3$ for $P_4$ and $b_2$ and $b_4$ for $B_5$. For $P_4$, the vertices $p_1$ and $p_4$ are called the \emph{leaf vertices}; and for $B_5$, the vertices $b_1$ and $b_5$ are the leaf vertices. Assume $G'$ to be any graph and $X$ to be a graph which is either a copy of $P_4$ or $B_5$. For a fixed vertex $q \in V(G')$, we define a new graph $G' \rhd_q X$ to be the graph obtained by identifying the vertex $q$ with an articulation vertex of $X$ (see Figures \ref{fig_G adjunct P4} and \ref{fig_G adjunct B5} for examples of $K_4 \rhd_q P_4$ and $K_4 \rhd_q B_5$, respectively). As a matter of reference, we call the new vertex in $G' \rhd_q X$ — obtained as a result of identifying two vertices — as the \emph{quotient vertex}; and continue to refer to the quotient vertex as $q$ itself. We now turn to our results. Firstly, it is easy to establish the following.

\begin{lemma} \label{lem_OLD_P4}
If $P$ is a $4$-path, then $\gamma^{OLD}(P)=4$.
\end{lemma}

\begin{lemma} \label{lem_OLD_B5}
If $B$ is a bull graph, then $\gamma^{OLD}(B)=3$.
\end{lemma}

\begin{proof}
Let $V(B)=\{ b_1, b_2, b_3, b_4, b_5 \}$, where $b_2$ and $b_4$ are the articulation vertices; and $b_1$ and $b_5$ are the leaf vertices of $B$. Then it is easy to check that $\{ b_2, b_3, b_4 \}$ is an OLD-code of $B$ and hence $\gamma^{OLD}(B) \leq 3$. See Figure \ref{fig_Bull_OLD} for the OLD-code demonstrated with black vertices in the figure.

\medskip
On the other hand, assume that $C$ is an OLD-code of $B$. Since $b_1$ and $b_5$ are degree $1$ vertices, their only neighbours, namely $b_2$ and $b_4$, respectively, must be in $C$ for the latter to be an open-dominating set of $B$. Moreover, at least one of $b_1$ and $b_3$ must be in $C$ for $b_2$ and $b_5$ to be open-separated in $B$. Hence, $|C| \geq 3$ and this establishes the result.
\end{proof}

This brings us to our result on the upper bound for OLD-numbers for block graphs.

\begin{theorem} \label{thm_OLD_ub}
Let $G$ be a connected open-twin-free block graph with no isolated vertices. Moreover, let $G$ neither be a copy of $P_2$ nor of $P_4$. Moreover, let $m_Q(G)$ be the number of non-leaf blocks of $G$ with at least one non-articulation vertex. Then $\gamma^{OLD}(G) \leq |V(G)|-m_Q(G) - 1$.
\end{theorem}

\begin{proof}
To start with, if $G$ is a copy of the bull graph, then $|V(G)|=5$ and $m_Q(G)=1$. Moreover, by Lemma \ref{lem_OLD_B5}, $\gamma^{OLD} (G)= 3 = |V(G)|- m_Q(G) -1$ and so, we are done. So, let us assume that $G$ is not a copy of the bull graph.

\medskip
We first choose a root block $K_0 \in \mathcal{K}(G)$ according to the following two possibilities.

\medskip
\textbf{Possibility 1}: $G \cong G' \rhd_q X$ for some block graph $G'$ and some graph $X$ that is a copy of either $P_4$ or $B_5$. In such a case, assume that $x$ is the articulation vertex of $X$ which is identified with the vertex $q$ of $G'$ to form $G$. Then, we choose $K_0$ to be the block of $G$ isomorphic to a $P_2$ with vertices $\{ x, z \}$, where $z$ is the leaf vertex of $X$ adjacent to $x$.

\medskip
\textbf{Possibility 2}: $G \not \cong G' \rhd_q X$ for all block graphs $G'$ and $X$ that is a copy of either $P_4$ or $B_5$. In this case, choose $K_0 \in \mathcal{K}(G)$ such that $|V(K_0)|= \min \{ |V(K)|: K \text{ is a leaf block of } G \}$.

Next, we construct a set $C \subset V(G)$ by the following rules.
\begin{enumerate}[{Rule} 1:]
	\item For every non-root leaf block $K \in \mathcal{K}(G)$, pick all vertices of $K$ in $C$.
	
	\item For every block $K \in \mathcal{K}(G)$ that is either the root block $K_0$ or is a non-leaf block in $\mathcal{K}(G) \setminus \{K_0 \}$, (i) pick all articulation vertices of $K$ in $C$; and (ii) pick all but one non-articulation vertices of $K$ in $C$.
\end{enumerate}
To compute the size of $C$, we note that, for the root block and every other non-leaf block $K$ with at least one non-articulation vertex, exactly one vertex is left out from it in $C$. This gives $|C|=|V(G)| - m_Q(G) - 1$. Thus, the upper bound for $\gamma^{OLD}(G)$ in the theorem is established on showing that $C$, indeed, is an OLD-code of $G$; which is what we prove next.

\medskip
To show that $C$ is an OLD-code of $G$, we notice first of all that, if the root block $K_0$ is isomorphic to $P_2$, then the (only) non-articulation vertex of $K_0$ is open-dominated by the articulation vertex of $K_0$ which belongs to $C$. In every other case when the root block is not isomorphic to $P_2$, all blocks $K \in \mathcal{K}(G)$ have $| V(K) \cap C | \geq 2$. This makes $C$ an open-dominating set of $G$. Now we show that $C$ is also an open-separating set of $V(G)$. So, let us assume that $u,v \in V(G)$ are distinct vertices of $G$. We now consider the following three cases.

\begin{case}[$u, v \in V(K)$ for some block $K \in \mathcal{K}(G)$]
We note here that, by the construction of $C$, for every block $K \in \mathcal{K}(G)$, at most one vertex of $K$ is not in $C$. This implies that at least one of $u$ and $v$, say $u$ without loss of generality, must be in $C$. Then $u$ open-separates $u$ and $v$ in $G$.
\end{case}

\begin{case}[$u \in V(K)$, $v \in V(K')$ for distinct $K,K' \in \mathcal{K}(G)$ and both $u,v \notin V(K) \cap V(K')$]
In this case, first let us assume that $V(K) \cap V(K') = \emptyset$. Now, both $K$ and $K'$ cannot be root blocks of $G$. Without loss of generality, let us assume that $K$ is not the root block of $G$. Since $|V(K) \cap C| \geq 2$ by the construction of $C$, the block $K$ has a vertex $x$ in $C$ other than $u$. Then, $x$ open-separates $u$ and $v$ in $G$. Therefore, assume that $V(K) \cap V(K') = \{ w \}$ for some vertex $w \in V(G)$. Now, if either one of $V(K)$ and $V(K')$, say $V(K)$ without loss of generality, has size at least $4$, then $|V(K) \cap C| \geq 3$ and so, at least one vertex in $V(K) \setminus \{ u, w \}$ belongs to $C$ which open-separates $u$ and $v$ in $G$. So, assume that $|V(K)| \leq 3$ and $|V(K')| \leq 3$. We next look at the following two subcases.

\begin{subcase}[$|V(K)|=|V(K')|=2$] Then, at least one of $K$ and $K'$ must be a non-leaf block, or else, $G$ would have open-twins. So, without loss of generality, suppose that $K$ is a non-leaf block. Then $\{ u \} = V(K) \cap V(K'')$ for some block $K'' ~(\neq K,K') \in \mathcal{K}(G)$. If however, $K'$ too is a non-leaf-block, then $\{ v \} = V(K') \cap V(K''')$ for some block $K''' ~(\neq K',K,K'') \in \mathcal{K}(G)$. Now, at least one of $K''$ and $K'''$ is not the root block. Without loss of generality, therefore, assume that $K''$ is not the root block. Then there is at least one vertex in $V(K'') \setminus \{ u \}$ which is in $C$ and, hence, open-separates $u$ and $v$ in $G$. So, let us assume that $K'$ is a leaf block, i.e. $v$ is a non-articulation vertex of $K'$.  Now again, if $K''$ is not the root block, then there is at least one vertex in $V(K'') \setminus \{ u \}$ which is in $C$ and, hence, open-separates $u$ and $v$ in $G$. So, now let us assume that $K''$ \emph{is} the root block. If $G$ satisfies the conditions of Possibility 1, then $|V(K'')| = 2$. However, if $G$ satisfies the conditions of Possibility 2, then again, since $K'$ is a leaf block and $|V(K')|=2$, by the minimality in size of the root block, we must have $|V(K'')|=2$. So, assume $z$ to be the non-articulation vertex of $K_0$.  If $P=G[z,u,w,v]$, we have $P \cong P_4$ with $u$ and $w$ being the articulation vertices and $z$ and $v$ being the leaf-vertices of $P$. Since $G \not \cong P_4$, we must have $G \cong G' \rhd_q P$ for some block graph $G'$ and some vertex $q \in \{ u,w \}$ (note that both $z$ and $v$ are non-articulation vertices of $G$). However, by the way we have chosen the root block $K_0$, we must have $q=u$. This implies that $u$ is the negative articulation vertex of some $K^* \in \mathcal{K}(G)$ such that $K^* \notin \{K, K_0 \}$. This, in turn, implies that $u$ and $v$ are open-separated in $G$ by some vertex in $(V(K^*) \cap C) \setminus \{ u \}$.
\end{subcase}

\begin{subcase}[At least one of $V(K)$ and $V(K')$ has size $3$] Without loss of generality, let us assume that $|V(K)| = 3$. So, assume that $K=G[w,u,y]$ for some vertex $y \in V(G)$. We must have $y\notin C$ (otherwise $y$ would open-separate $u$ and $v$).
We first assume that $K$ is the root block. If $v$ is an articulation vertex of $K'$, then $\{ v \} = V(K') \cap V(K''')$ for some block $K''' ~(\neq K') \in \mathcal{K}(G)$. This implies that there exists a vertex in $V(K''') \setminus \{ v\}$ which open-separates $u$ and $v$ in $G$. Moreover, if $v$ is a non-articulation vertex of $K'$, then we must have $|V(K')|=3$ (or else, $K'$ is a leaf block of size smaller than the root block which is a contradiction). Assume that $V(K')=\{ w,v,a \}$ for some vertex $a \in V(G)$. If $a$ is also a non-articulation vertex of $K'$, then $K'$ is a leaf block that is not a root block and, hence, $a \in C$. If however, $a$ is an articulation vertex of $K'$, then also, $a \in C$. Thus, either way, $a$ open-separates $u$ from $v$ in $G$. Thus, we are done in the case that $K$ is the root block of $G$.

\medskip
So, let us now assume that $K$ is \emph{not} the root block of $G$. Now, if $y \in C$, then $y$ open-separates $u$ and $v$ in $G$. So, let us assume that $y \notin C$, which implies that $y$ is a non-articulation vertex of $K$. This, in turn, implies that $u$ is an articulation vertex of $K$, or else, $K$ would be a leaf block of $G$ that is not the root block and so, $y \in C$, contrary to our assumption. So, let $\{ u \} = V(K) \cap V(K'')$ for some block $K'' ~(\neq K) \in \mathcal{K}(G)$. If however, $v$ is also an articulation vertex of $K'$, then $\{ v \} = V(K') \cap V(K''')$ for some block $K''' ~(\neq K') \in \mathcal{K}(G)$. Now, at least one of $K''$ and $K'''$ is not the root block. So, without loss of generality, assume that $K''$ is not the root block. Then, there is at least one vertex in $V(K'') \setminus \{ u \}$ which is in $C$ and, hence, open-separates $u$ and $v$ in $G$. So, let us assume that $v$ is a non-articulation vertex of $K'$. Again, if $K''$ is not the root block, then there exists at least one vertex in $V(K'') \setminus \{ u \}$ which is in $C$ and, hence, open-separates $u$ and $v$ in $G$. So, now assume $K''$ to \emph{be} the root block. If $|V(K'')| \geq 3$, then there exists a vertex of $V(K'') \setminus \{ u \}$ in $C$ and, hence, open-separates $u$ and $v$ in $G$. So, let us assume that $|V(K'')|=2$ and that $z$ is the non-articulation  of $K''$. If $|V(K')|=3$, then suppose that $V(K')=\{ w,v,a \}$ for some vertex $a \in V(G)$. If $a$ is also a non-articulation vertex of $K'$, then $K'$ is a leaf block that is not a root block and hence, $a \in C$. If however, $a$ is an articulation vertex of $K'$, then also, $a \in C$. Thus, either way, $a$ open-separates $u$ from $v$ in $G$ and we are done in the case that $|V(K')|=3$. So, let us finally assume that $|V(K')|=2$ and that $v \in V(K')$ is a non-articulation vertex of $K'$.

\medskip
If $B=G[z,u,y,w,v]$, we have $B \cong B_5$ with $u$ and $w$ being the articulation vertices and $z$ and $v$ being the leaf-vertices of $B$. Since by our assumption, $G \not \cong B_5$, we have $G \cong G' \rhd_q B$ for some block graph $G'$ and some vertex $q \in \{ u,w \}$ (note that $z$, $y$ and $v$ are non-articulation vertices of $G$). Now, by the way we have chosen the root block, this implies that we must have $q=u$. This further implies that $u$ is the negative articulation vertex of $K^*$ for some $K^* \in \mathcal{K}(G)$ such that $K^* \notin \{K, K_0 \}$. Hence, $u$ and $v$ are open-separated in $G$ by some vertex in $V(K^*) \cap C \setminus \{ u \}$.
\end{subcase}
\end{case}

This, therefore, proves that $C$ is an OLD-code of $G$ and with that, we prove the theorem.
\end{proof}

Applying Theorem \ref{thm_OLD_ub} to each connected component of a block graph, one has the following general result.

\begin{corollary} \label{cor_OLD_ub}
Let $G$ be an open-twin-free block graph with $k$ connected components and no isolated vertices. Moreover, let no component of $G$ be either a copy of $P_2$ or of $P_4$. Also, let $m_Q(G)$ be the number of non-leaf blocks of $G$ with at least one non-articulation vertex. Then $\gamma^{OLD}(G) \leq |V(G)|-m_Q(G) - k$.
\end{corollary}

Foucaud et al. \cite{foucaud2021characterizing} have shown that, for any open-twin-free graph $G$, $\gamma^{OLD} (G) \leq |V(G)| - 1$ unless $G$ is a special kind of bipartite graph called half-graph (a \emph{half-graph} is a bipartite graph with both parts of the same size, where each part can be ordered so that the open neighbourhoods of consecutive vertices differ by exactly one vertex~\cite{EH84}). Noting that $P_2$ and $P_4$ are the only block graphs that are half-graphs, Theorem~\ref{thm_OLD_ub} can be seen as a refinement of this result for block graphs.

\begin{figure}[h!]
\vspace{-2mm}
     \centering
         \includegraphics[scale=0.42]{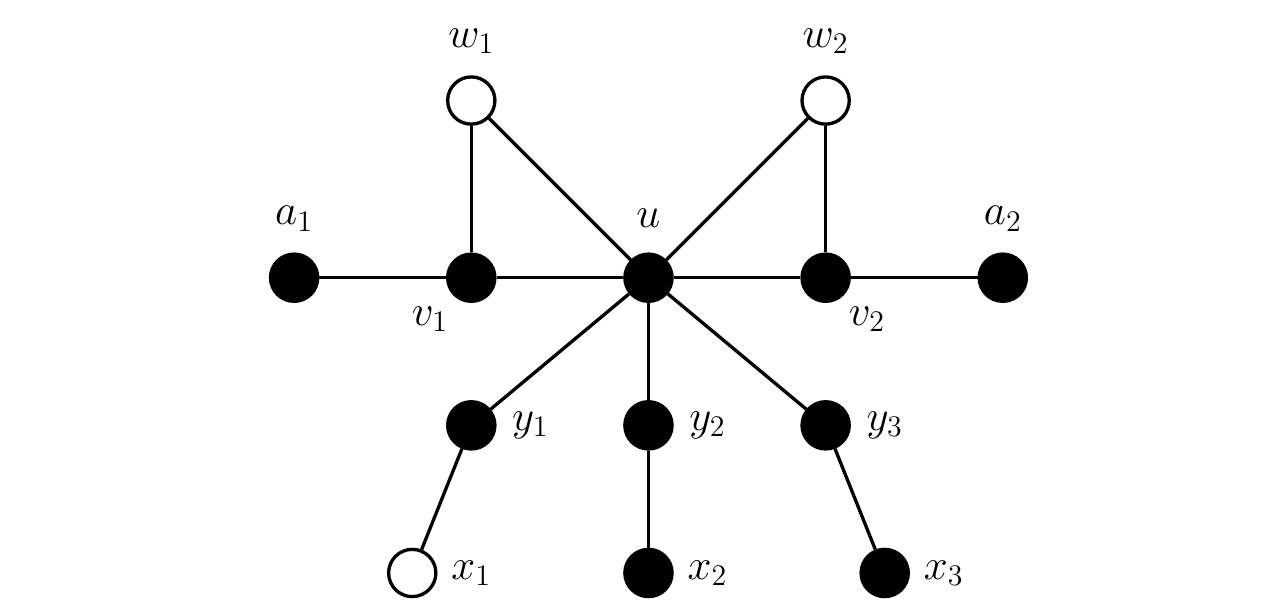}
     \caption{Graph $G_{2,3}$ whose OLD-number attains the upper bound in Theorem \ref{thm_OLD_ub}. The black vertices represent those included in the OLD-code $C$ of $G$ as described in the proof of the Theorem \ref{thm_OLD_ub}.}
     \label{fig_OLD_ub}\vspace*{-2mm}
\end{figure}

\medskip
We now show that the upper bound given in Theorem \ref{thm_OLD_ub} is tight and is attained by arbitrarily large connected block graphs.
To prove so, for two non-negative integers $k$ and $l$ such that $k+l \geq 2$, let us define a subclass $G_{k,l}$ of block graphs by the following rule: Let $T_1, T_2, \ldots, T_k$ be $k$ copies of $K_3$ with $V(T_i)=\{ u_i, v_i, w_i \}$ for each $1 \leq i \leq k$. Further, let $A_1, A_2, \ldots , A_k$ are $k$ copies of $P_2$ with $V(A_i)=\{ a_i, b_i \}$ for each $1 \leq i \leq k$, and $L_1, L_2, \ldots , L_l$ be $l$ copies of $P_3$ with $V(L_j)=\{ x_j, y_j, z_j \}$ for each $1 \leq j \leq l$. Let $G_{k,l}$ be the graph obtained by identifying the vertices $v_i$ with $b_i$ for each $1 \leq i \leq k$ and identifying the vertices $u_i$ and $z_j$ for all $1 \leq i \leq k$ and $1 \leq j \leq l$ into a single vertex $u$, say. See Figure \ref{fig_OLD_ub} for an example of a construction of $G_{k,l}$ with $k=2$ and $l=3$. As a matter of reference, we continue to call the vertices of $G_{k,l}$ obtained by identifying $v_i$ with $b_i$, for all $1 \leq i \leq k$, as $v_i$ itself.

\begin{proposition}
For all positive integers $k$ and $l$ with $l+k \geq 2$, $\gamma^{OLD} (G_{k,l})=|V(G_{k,l})| - m_Q(G_{k,l}) - 1$.
\end{proposition}

\begin{proof}
First of all, we have $|V(G_{k,l})| = 3k + 2l + 1$, and $m_Q(G_{k,l}) = k$. Therefore, by Theorem \ref{thm_OLD_ub}, we have $\gamma^{OLD} (G_{k,l}) \leq 2(k+l)$.

\medskip
Now, let $C$ be an OLD-code of $G_{k,l}$. For $C$ to be an open-dominating set of $G_{k,l}$, the only neighbours $v_i$ and $y_j$ of the degree $1$ vertices $a_i$ and $x_j$, respectively, of $G_{k,l}$ must be in $C$. Moreover, $u \in C$ for each pair of $w_i$ and $a_i$ to be open-separated in $G_{k,l}$ by $C$. Similarly, at least all but one of the $x_i$'s must belong to $C$ for each pair of vertices $y_i$ and $y_j$, for $1 \leq i < j \leq l$, to be open-separated in $G_{k,l}$ by $C$. Let us first assume that $x_2, x_3, \ldots , l_l$ belong to $C$, without loss of generality. Then, for all $1 \leq i \leq k$, at least one of $w_i$ and $a_i$ must be in $C$ for the pair $y_1$ and $v_i$, for $1 \leq i < j \leq k$, to be open-separated in $G_{k,l}$ by $C$. Adding up, therefore, we have $|C| \geq 2(k+l)$. On the other hand, if $x_1, x_2, \ldots , x_l \in C$, then, to open-separate the vertices $v_i$ and $v_j$ (for $1 \leq i < j \leq k$), we need at least one of $w_i$ and $a_i$ (for $k-1$ such $i$) to belong to $C$. Then again, the count adds up to $|C| \geq 2(k+l)$. This proves the proposition.
\end{proof}

\section{Lower bounds} \label{sec:LB}

The general lower bound for the size of an identifying code using the number of vertices is $\gamma^{ID}(G)\geq \lceil \log_2(|V(G)|+1)\rceil$ \cite{karpovsky1998new}.
However, to reach this bound, a graph needs to have a large VC-dimension~\cite{bousquet2015identifying} (the {\em VC-dimension} of a graph $G$ is the size of a largest \emph{shattered set}, that is, a set $S$ of vertices such that for every subset $S$' of $S$, some closed neighbourhood in $G$ intersects $S$ exactly at $S'$). Indeed, if a graph has VC-dimension $c$, then any identifying code has size at least $O(|V(G)|^{1/c})$~\cite{bousquet2015identifying}. The value $1/c$ is not always tight for specific subclasses of graphs of VC-dimension $c$, see for example the case of line graphs, which have VC-dimension at most~$4$ but for which the tight order for the lower bound is $\Omega(|V(G)|^{1/2})$~\cite{foucaud2013identifying}. Similar results hold for LD- and OLD-codes, by using the same techniques as in~\cite{bousquet2015identifying}.

\medskip
Block graphs have VC-dimension at most~$2$ (one can check that a shattered set of size~3 would imply the existence of an induced 4-cycle or diamond), and thus, using the result from~\cite{bousquet2015identifying}, their ID-number is lower bounded by $\Omega(|V(G)|^{1/2})$. In this section, we improve this lower bound to a linear one, and give a tight result.

For the rest of this article, given a block graph $G$, by $\mathcal{K}_{leaf}(G)$ we shall mean the set of all leaf blocks of $G$ with at least one edge in the block. Moreover, by the symbol $n_i(G)$, we shall mean the number of vertices of degree $i$ in the graph $G$.

\begin{lemma}\label{newlem_block}
Let $G$ be a connected block graph with at least one edge and with blocks $B_1, B_2, \ldots B_h$, say. Then, there exist distinct vertices $v_0, v_1, v_2, \ldots , v_h$ of $G$ such that $v_0, v_1 \in V(B_1)$ and $v_i \in V(B_i)$ for all $2 \leq i \leq h$.
\end{lemma}

\begin{proof}
Since $G$ has at least one edge, note that every block of $G$ has at least one edge. The proof is by induction on $h$ with the base case being $h=1$. In the base case, we have $G=B_1$. Since $G$ has at least one edge $v_0v_1$, say, the two distinct vertices are $v_0$ and $v_1$ and we are done. So, let us assume the induction hypothesis for any connected block graph $G'$ with at least one edge and with $h'$ blocks, where $1 \leq h' \leq h-1$. Therefore, let $h \geq 2$. Without loss of generality, let us assume that $B_h$ is a leaf block of $G$. Let $G'=G-( V(B_h) \setminus art^-(B_h) )$. Then $G'$ is also a connected block graph on $h-1$ blocks. Since $G'$ contains the block $B_1$, the former has at least one edge. Then by the induction hypothesis, there exist distinct vertices $v_0, v_1, v_2, \ldots , v_{h-1}$ of $G$ such that $v_0, v_1 \in V(B_1)$ and $v_i \in V(B_i)$ for all $2 \leq i \leq h-1$.  Since $B_h$ has at least one edge, it implies that there exists a vertex $v_h \in V(B_h) \setminus art^-(B_h)$. This implies that $v_h \neq v_i$ for all $1 \leq i \leq h-1$ and hence, the result holds.
\end{proof}

\begin{lemma}\label{lem_K(G) leq V(G)-1}
Let $G$ be a connected block graph with at least one edge. Then we have
$$|\mathcal{K} (G)| \leq |V(G)| - 1 - |\mathcal{K}_{leaf} (G)| + n_1(G).$$
\end{lemma}

\begin{proof}
Let $\mathcal{L}(G) = \{ L \in \mathcal{K}_{leaf} (G) : L \cong K_2 \}$ and  $G^*$ be a graph obtained from $G$ by, for each $L \in \mathcal{L}(G)$, introducing a new vertex and making it adjacent to both elements of $V(L)$. Thus, $G^*$ is a block graph in which every leaf block has at least $3$ vertices. We also note here that
\begin{enumerate}[(1)]
\itemsep=0.9pt
	\item $|\mathcal{L} (G)| = n_1 (G)$,
	
	\item $|V(G^*)| = |V(G)|+|\mathcal{L} (G)| = |V(G)| + n_1 (G)$,
	
	\item $|\mathcal{K} (G)| = |\mathcal{K} (G^*)|$ and that
	
	\item $|\mathcal{K}_{leaf} (G)| = |\mathcal{K}_{leaf} (G^*)|$.
\end{enumerate}

Now, let $|\mathcal{K} (G^*)| = h$ and $\mathcal{K} (G^*) = \{B_1, B_2, \ldots , B_h \}$. Then, by Lemma \ref{newlem_block}, there exist distinct vertices $v_0, v_1, v_2, \ldots , v_h$ of $G$ such that $v_0, v_1 \in V(B_1)$ and $v_i \in V(B_i)$ for all $2 \leq i \leq h$. Moreover, since at least one vertex in each leaf block of $G^*$ is not any of the vertices $v_i$, we have
\begin{flalign*}
|\mathcal{K} (G)| = |\mathcal{K} (G^*)| = h = |\{ v_0, v_1, v_3, \ldots , v_h \}| - 1 &\leq |V(G^*)| - |\mathcal{K}_{leaf} (G^*)| - 1\\
&= |V(G)|+ n_1(G) - |\mathcal{K}_{leaf} (G)| - 1.
\end{flalign*}

\vspace*{-6mm}
\end{proof}

\begin{corollary} \label{cor_K(G)-n_0(G)}
Let $G$ be a block graph with $k$ connected components. Then, we have
$$|\mathcal{K} (G)| - n_0(G) \leq |V(G)| - k - |\mathcal{K}_{leaf} (G)| + n_1(G).$$
\end{corollary}

\begin{proof}
Assume that $k=p+q$ such that $G_1, G_2, \ldots , G_p$ are the connected components of $G$, each with at least one edge; and that $S_1, S_2, \ldots , S_q$ are the components of $G$, each with a single vertex. Then, we have
\begin{flalign*}
|\mathcal{K}(G)| &= \sum_{1 \leq i \leq p} |\mathcal{K}(G_i)| + \sum_{1 \leq j \leq q} |\mathcal{K}(S_j)|\\
&\leq q-p + \sum_{1 \leq i \leq p} \Big(|V(G_i)| - |\mathcal{K}_{leaf}(G_i)| + n_1(G_i) \Big) &\text{[using Lemma \ref{lem_K(G) leq V(G)-1}]\hspace*{10mm}}\\
&= |V(G)|- k - |\mathcal{K}_{leaf}(G)| + n_1(G) + n_0(G) &\text{[since $q=n_0(G)$]\hspace*{10mm}}
\end{flalign*}

\vspace*{-7mm}
\end{proof}

\begin{corollary} \label{cor_K(G) <= V-k}
Let $G$ be a block graph with $k$ connected components. Then we have
$$|\mathcal{K} (G)| - n_0(G) \leq |V(G)| - k.$$
\end{corollary}

\begin{proof} \label{cor_K(G) leq V(G)-1}
The result follows from Corollary \ref{cor_K(G)-n_0(G)} and the fact that $n_1(G) \leq |\mathcal{K}_{leaf} (G)|$.
\end{proof}

Before we come to our results, we define the following notations. For a given code $C$ of a connected block graph $G$, let us assume that the sets $C_1,C_2,...,C_k$ partition the code $C$ such that the induced subgraphs $G[C_1], G[C_2], \ldots , G[C_k]$ of $G$ are the $k$ connected components of the subgraph $G[C]$ of $G$ induced by $C$. Note that each $C_i$ is a block graph (since every induced subgraph of a block graph is also a block graph).

\begin{definition}\label{def_V_i}
Given a block graph $G$, its vertex set $V(G)$ is partitioned into the four following parts.
\begin{enumerate}[(1)]
\itemsep=0.9pt
\item $V_1=C$,
\item $V_2 = \{ v \in V(G) \setminus V_1 : |N_G(v) \cap C| = 1 \}$,
\item $V_3 = \{ v \in V(G) \setminus V_1 : \text{ there exist distinct } i, j \leq k \text{ such that } N_G(v) \cap C_i \neq \emptyset \text{ and } \\
          \hspace*{10mm}  N_G(v) \cap C_j \neq \emptyset \}$,
\item $V_4 = V(G) \setminus (V_1 \cup V_2 \cup V_3)$. Note that, for all $v \in V_4$, we have $N_G(v) \cap C \subset C_i$ for some $i$ and that $|N_G(v) \cap C_i| \geq 2$.
\end{enumerate}
\end{definition}

We now prove a series of lemmas establishing upper bounds on the sizes of each of the vertex subsets $V_1, V_2, V_3$ and $V_4$ of a connected block graph $G$.

\begin{lemma}\label{claimV2}
Let $G$ be a connected block graph and $C$ be a code of $G$. Then the following are upper bounds on the size of the vertex subset $V_2$ of $G$.
  \begin{enumerate}[(1)]
  \itemsep=0.9pt
    \item $|V_2| \leq |C|-n_0(G[C])$ if $C$ is an ID-code.

    \item $|V_2| \leq |C|$ if $C$ is an LD-code.

    \item $|V_2| \leq |C|-n_1(G[C])$ if $C$ is an OLD-code.
\end{enumerate}
\end{lemma}

\begin{proof}
By definition of $V_2$, each vertex $v \in V_2$ has a unique neighbor $u$ in $C$, i.e. $N_G(v) \cap C = \{u\}$. Hence, there can be at most $|C|$ vertices in $V_2$ and this proves (2).

\smallskip
If $C$ is an ID-code, $u$ cannot be isolated in $G[C]$ (or else, $u$ and $v$ will not be closed-separated in $G$ by $C$). Thus, there are at most $|C|-n_0(G[C])$ vertices in $V_2$ and this proves (1).

\medskip
Finally, if $C$ is an OLD-code and $u$ has a neighbor $w \in N_G(u) \cap C$ such that $deg_{G[C]} (w) = 1$, then $v$ and $w$ are not open-separated. Thus, there  are at most $|C|-n_1(G[C])$ vertices in $V_2$ and this proves (3).
\end{proof}

\begin{definition} \label{def_aux}
Given a connected block graph $G$ and a code $C$ of $G$, let $G[C_1], G[C_2], \ldots, G[C_k]$ be all the connected components of $G[C]$, where $C_1, C_2, \ldots , C_k$ are subsets of $C$. Moreover, let the vertex set $V(G)$ be partitioned into the subsets $V_1, V_2, V_3, V_4$ as in Definition \ref{def_V_i}. Then, consider the bipartite graph $F_C(G)$, where $A = \{a_j: v_j \in V_3 \}$ and $B =  \{u_i: G[C_i] \text{ is a connected component of }\\ G[C] \}$ are the two \emph{parts} of $V(F_C(G))$. As for the edge set $E(F_C(G))$, for each vertex $v_j$ in $V_3$, we add an edge between $a_j$ and $u_i$ if $v_j$ is adjacent to a vertex in $C_i$. The graph $F_C(G)$ is called the \emph{auxiliary graph} of $G$.
\end{definition}

\begin{lemma} \label{lem_F_C(G)}
For a connected block graph $G$ and a code $C$ of $G$, the auxiliary graph $F_C(G)$ is a forest.
\end{lemma}

\begin{proof}
If there is a cycle in $F_C(G)$, there would be a cycle in $G$ involving two vertices of different connected components $G[C_i]$ and $G[C_j]$, say. By the definition of a block graph, the latter cycle in $G$ has to induce a complete subgraph in $G$. However, that would imply that $G[C_i]$ and $G[C_j]$ must be the same component of $G[C]$ which is a contradiction. Thus, $F_C(G)$ is cycle-free and, hence, is a forest.
\end{proof}

\begin{lemma} \label{claimV3}
Let $G$ be a connected block graph, $C$ be a code of $G$ and $G[C_1], G[C_2], \ldots, G[C_k]$ be all the connected components of $G[C]$, where $C_1, C_2, \ldots , C_k$ are subsets of $C$. Then, we have $|V_3| \leq k-1$.
\end{lemma}

\begin{proof}
By Lemma \ref{lem_F_C(G)}, $F_C(G)$ is a forest. Let $V(F_C(G)) = A \sqcup B$ be as defined above. Then we have $|B|=k$. Therefore, $F_C(G)$, on account of being a forest, has at most $|A|+k-1$ edges. Also, since $F_C(G)$ is bipartite, its number of edges is $\sum_{a \in A} deg_{F_C(G)}(a) \geq 2|A|$ (the last inequality holds since any vertex in the part $A$ of $V(F_C(G))$ is adjacent to at least two distinct vertices of $B$). Thus, the result follows from the fact that $|V_3| = |A|$.
\end{proof}

\begin{lemma}\label{claimV4}
Let $G$ be a connected block graph, $C$ be a code of $G$ and $G[C_1], G[C_2], \ldots, G[C_k]$ be all the connected components of $G[C]$, where $C_1, C_2, \ldots , C_k$ are subsets of $C$. Then, we have $|V_4| \leq |C|-k$. In particular,
  \begin{enumerate}[(1)]
\item $|V_4| \leq |\mathcal{K}(G[C])| - |\mathcal{K}_{leaf}(G[C])| \leq |C|-3k$ if $C$ is an ID-code;

\item $|V_4| \leq |\mathcal{K}(G[C])| \leq |C|- 2k +n_1(G[C])$ if $C$ is an OLD-code.
  \end{enumerate}
\end{lemma}

\begin{proof}
Let $v$ be any vertex in $V_4$ and let $G[C_i]$ be the component of $G[C]$ such that $N_G(v)\cap C\subseteq C_i$. Moreover, $|N_G(v) \cap C_i| \geq 2$. Then, notice that $N_G(v) \cap C$ must be a subset of exactly one block of $G[C_i]$, or else, $G[C_i]$ would be disconnected, as $v \notin C$. This implies that $|V_4| \leq |\mathcal{K}(G[C])| - n_0(G[C]) \leq |C| - k$, by Corollary \ref{cor_K(G) <= V-k}. We now prove the more specific bounds for ID- and OLD-codes.

\medskip
(1) First consider the case where $C$ is an ID-code. Let $G[C_i]$ be a connected component of $G[C]$ such that at least one vertex of $V_4$ is adjacent to some vertices in $C_i$. In particular, $|C_i| \geq 2$ and since $C$ is a ID-code, $G[C_i]$ is closed-twin-free. We first show the following.

\begin{claim}\label{subclaim_V4.1}
No element of $V_4$ is adjacent to the vertices of the leaf blocks of $G[C_i]$.
\end{claim}

\begin{proofofclaim}
Suppose that $L$ is a leaf block of $G[C_i]$. Then we must have $L \cong K_2$, or else, at least two vertices in $V(L)$ are not closed-separated in $G$ by $C$. So, assume that $V(L) = \{ x,y \}$. Then at least one of $x$ and $y$ must be a non-articulation vertex of $G[C_i]$. Without loss of generality, suppose that $y$ is a non-articulation vertex of $G[C_i]$. If there exists a vertex $v$ of $V_4$ such that $N_G(v)\cap C=V(L)$, then $v$ and $y$ would not be closed-separated in $G$ by $C$ which is a contradiction. Hence, no element of $V_4$ is adjacent to the vertices of the leaf blocks of $G[C_i]$.
\end{proofofclaim}

This implies that the number of vertices of $V_4$ that can be adjacent to the vertices of $G[C_i]$ are at most $|\mathcal{K}(G[C_i])| - |\mathcal{K}_{leaf} (G[C_i])|$. Now, we must have the following.

\begin{claim}\label{subclaim_V4.2}
$|\mathcal{K}_{leaf} (G[C_i])| \geq 2$
\end{claim}

\begin{proofofclaim}
If, on the contrary, $|\mathcal{K}_{leaf}(G[C_i])|=1$, then $|C_i|=1$; or else, all pairs of vertices of $G[C_i]$ are not closed-separated in $G$ by $C$. This contradicts the fact that $|C_i| \geq 2$.
\end{proofofclaim}

Therefore, by the above two claims and the fact that any vertex of $V_4$ having its neighbors in a component $G[C_i]$ of $G[C]$ is adjacent to the vertices of exactly one block of $G[C_i]$, the number of vertices of $V_4$ adjacent to the vertices of $G[C_i]$ is at most $|\mathcal{K}(G[C_i])| - |\mathcal{K}_{leaf} (G[C_i])| \leq |\mathcal{K} (G[C_i])| - 2 \leq |C_i| - 3$ (the last inequality is by the fact that $G[C_i]$ is connected, that is, $k=1$ and $n_0(C_i) = 0$ in Corollary \ref{cor_K(G) <= V-k}). Hence,
$$|V_4| \leq \sum_{1 \leq i \leq k} \Big( |\mathcal{K}(G[C_i])| - |\mathcal{K}_{leaf} (G[C_i])| \Big) \leq \sum_{1 \leq i \leq k} \Big( |C_i| - 3 \Big) =  |C| - 3k.$$

(2) In the case that $C$ is an OLD-code of $G$, we have $n_0(G[C]) = 0$. So, assume that $G[C_i]$ is a connected component of $G[C]$ with at least one edge. Then, by Lemma \ref{lem_K(G) leq V(G)-1}, we have $|\mathcal{K} (G[C_i])| \leq |C_i| - 1 - |\mathcal{K}_{leaf} (G[C_i])| + n_1(G[C_i]) \leq |C_i| - 2 + n_1(G[C_i])$. This implies that 
$$|V_4| \leq |\mathcal{K}(G[C])| = \sum_{1 \leq i \leq k} |\mathcal{K}(G[C_i])| \leq \sum_{1 \leq i \leq k} \Big(|C_i| - 2 + n_1(G[C_i]) \Big) = |C|- 2k + n_1(G[C]).$$

This proves the lemma.
\end{proof}

\begin{theorem}\label{thm_lb}
Let $G$ be a connected block graph. Then we have
\begin{itemize}
   \item  $\gamma^{ID}(G) \geq \frac{|V(G)|}{3}+1$;

   \item  $\gamma^{LD}(G) \geq \frac{|V(G)|+1}{3}$; and

   \item $\gamma^{OLD}(G) \geq \frac{|V(G)|+2}{3}$.
\end{itemize}
\end{theorem}

\begin{proof}
Let $|V(G)| = n$. Assume $C$ to be a code of $G$ and $G[C_1], G[C_2], \ldots, G[C_k]$ be all the connected components of $G[C]$, where $C_1, C_2, \ldots , C_k$ are subsets of $C$. Recalling from Definition \ref{def_V_i} the sets $V_1, V_2, V_3, V_4$ that partition $V(G)$, we prove the theorem using the relation $|V(G)|=|C|+|V_2|+|V_3|+|V_4|$ and the upper bounds for $|V_2|$, $|V_3|$ and $|V_4|$ in Lemmas \ref{claimV2}, \ref{claimV3} and \ref{claimV4}, respectively.

\medskip
If $C$ is an ID-code, then we have
$$
\begin{array}{rcl}
n & =    & |C|+|V_2|+|V_3|+|V_4| \\
  & \leq & |C|+|C|-n_0(G[C])+k-1+|C| - 3k\\
  & \leq    & 3|C|-2k-1 \leq 3|C|-3, \quad \text{using $k \geq 1$}.
\end{array}
$$
Hence, the result holds.
\eject

\noindent If $C$ is an LD-code, then
$$
\begin{array}{rcl}
n & =    & |C|+|V_2|+|V_3|+|V_4| \\
  & \leq & |C|+|C|+k-1+|C|- k \\ 
  & =    & 3|C|- 1 \\ 
\end{array}
$$
and, hence, the result holds.

\medskip
Finally, if $C$ is an OLD-code, then we have
$$
\begin{array}{rcl}
n & =    & |C|+|V_2|+|V_3|+|V_4| \\
  & \leq & |C|+|C|-n_{1}(G[C])+k-1+|C| -2k +n_1(G[C])\\
  & =    & 3|C|- k -1 \leq 3|C|-2, \quad \text{using $k \geq 1$}.
\end{array}
$$

Hence, the result holds.
\end{proof}

\begin{figure}[b!]
     \centering
     \begin{subfigure}[t]{0.24\textwidth}
         \centering
         \includegraphics[angle=90, scale=1.4]{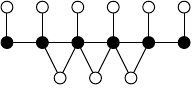}
         \caption{\footnotesize}
         \label{fig_ID_extremal}
     \end{subfigure}
     \begin{subfigure}[t]{0.24\textwidth}
         \centering
         \includegraphics[angle=90, scale=1.4]{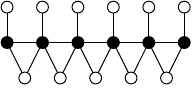}
         \caption{\footnotesize}
         \label{fig_LD_extremal}
     \end{subfigure}
    \begin{subfigure}[t]{0.24\textwidth}
         \centering
         \includegraphics[scale=0.35]{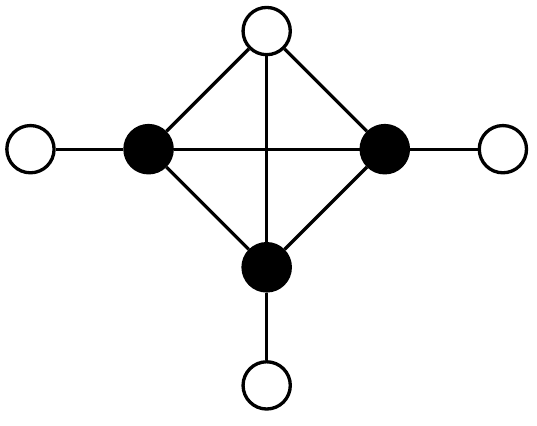}
         \caption{\footnotesize}
         \label{fig_OLD_extremal}
     \end{subfigure}
     \begin{subfigure}[t]{0.24\textwidth}
         \centering
         \includegraphics[angle=90, scale=1.4]{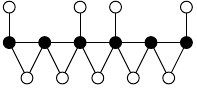}
         \caption{\footnotesize}
         \label{fig_OLD2_extremal}
     \end{subfigure}\vspace*{-2mm}
     \caption{Extremal cases where the lower bounds in Theorem \ref{thm_lb} are attained. The black vertices form a minimum (a) ID-code, (b) LD-code, (c) and (d) OLD-code.}
     \label{extremal}
\end{figure}

\medskip
We now look at examples of which connected block graphs are extremal with respect to the bounds of the previous theorem. To that end, we show that there are infinitely many connected block graphs whose ID- and LD-numbers reach the first and the second lower bounds, respectively, of Theorem~\ref{thm_lb}. In contrast however, there is exactly one connected block graph whose OLD-number reaches the last lower bound in Theorem~\ref{thm_lb}. In order to calculate the code numbers in the extremal cases attaining the bounds of Theorem~\ref{thm_lb}, the inequalities in the previous calculations in some of the earlier results leading to the bounds will have to be changed to equalities.

\begin{proposition} \label{prop_ID_extremal}
For $k \geq 3$, consider a path on vertices $u_1, u_2, \ldots,u_\ell$. For all $1 \leq i \leq \ell$, attach a vertex $v_i$ by the edge $v_iu_i$ and, for each pair $u_i,u_{i+1}$ for $\ell \geq 4$ and $2 \leq i \leq \ell-2$, attach a vertex $w_i$ by edges $w_iu_i$ and $w_iu_{i+1}$. See Figure \ref{fig_ID_extremal} with $\ell=6$. Then, such graphs are the only ones whose ID-numbers attain the lower bound in Theorem \ref{thm_lb}.
\end{proposition}

\begin{proof}
For extremal graphs whose ID-numbers attain the lower bound in Theorem \ref{thm_lb}, we have equalities in the equations in the proof of Theorem \ref{thm_lb} when $C$ is an ID-code. We therefore have in this case $k=1$ (that is, $G[C]$ is connected) which implies that $|V_2| = |C|$, $|V_3|=0$ and $|V_4|=|C|-3$. Tracing back the equality for $|V_4|$, it stems from equality in the statement of Lemma \ref{claimV4}(1), that is,
$$|V_4| = |\mathcal{K}(G[C])| - |\mathcal{K}_{leaf}(G[C])| = |C|-3k.$$

This further implies $|\mathcal{K}_{leaf}(G[C])| = 2$ in Claim B of the proof of Lemma \ref{claimV4}. Recall that $k=1$. Hence, putting these values in the preceding equation, we have $|V_4| = |\mathcal{K}(G[C])| - 2 = |C|-3$ or $|\mathcal{K}(G[C])| = |C|-1$. This implies that $G[C]$ must be a tree and more particularly, a path, since $|\mathcal{K}_{leaf}(G[C])| = 2$. Moreover, with the values of $|V_2| = |C|$ and $|V_4| = |C|-3$ coupled with the fact that no vertex of $V_4$ can have neighbors in leaf blocks of $C$ (Claim A in proof of Lemma \ref{claimV4}), any extremal graph with respect to the lower bound for ID-numbers in Theorem \ref{thm_lb} must be as described in the statement of the proposition.
\end{proof}

\begin{proposition} \label{prop_LD_extremal}
There exist arbitrarily large connected block graphs whose LD-number attains the lower bound in Theorem \ref{thm_lb}.
\end{proposition}

\begin{proof}
For any $\ell \geq 1$, consider a path on vertices $u_1, u_2, \ldots,u_\ell$. For all $1 \leq i \leq \ell$, attach a vertex $v_i$ by the edge $v_iu_i$ and, for each pair $u_i,u_{i+1}$ for $\ell \geq 2$ and $1 \leq i \leq \ell-1$, attach a vertex $w_i$ by edges $w_iu_i$ and $w_iu_{i+1}$. We call the graph $G$. See Figure \ref{fig_LD_extremal} for an example with $k=6$. Then it can be verified that $C = \{u_1, u_2, \ldots , u_\ell\}$ is a minimum LD-code of $G$. Since $|V(G)| = 3\ell-1$, the LD-number of $G$ attains the lower bound.
\end{proof}

\begin{proposition} \label{prop_OLD_extremal}
(1) Consider the block graph $Z$ consisting of a clique on four vertices, of which three are support vertices to one leaf each. See Figure \ref{fig_OLD_extremal}. Then $Z$ is the only block graph whose OLD-number attains the bound in Theorem \ref{thm_lb}.

\medskip
\noindent (2) For any $G \not \cong Z$, the bound for the OLD-number in Theorem \ref{thm_lb} becomes
$$\gamma^{OLD} (G) \geq \frac{|V(G)|}{3}+1.$$
In this case, there are arbitrarily large connected block graphs whose OLD-number attains this lower bound.
\end{proposition}

\begin{proof}
For any block graph $G$ whose OLD-number attains the bound in Theorem \ref{thm_lb}, we must have equalities in the equations in the proof of Theorem \ref{thm_lb} when $C$ is an OLD-code. This implies that we must have $k=1$. Moreover, we must have $|V_2| = |C| - n_1(G[C])$, $|V_3|=0$ and
\begin{flalign}
|V_4| = |C| - 2 + n_1(G[C]) = |\mathcal{K}(G[C])|
\qquad \text{(using Lemma \ref{claimV4} and $k=1$)} \nonumber
\end{flalign}

\noindent (1) The last equation stems from an equality in the statement of Lemma \ref{lem_K(G) leq V(G)-1}, which implies that $|\mathcal{K}_{leaf}(G[C])|=1$. This implies that $C$ is a clique. Since $C$ is an OLD-code, we must have $|C| \geq 2$. If $|C| = 2$, then we have $|V_4| = n_1(G[C]) = 1$. This is a contradiction, as $|C| = 2$ implies that $G[C] \cong P_2$ and hence, $n_1(G[C]) = 2$. Thus, we  have $|C| \geq 3$. This implies that $n_1(G[C]) = 0$ and so, we get $|V_4| = |C|-2 = 1 \implies |C|=3$. Moreover, $|V_2| = 3$ determines the graph to be $Z$.

\medskip
\noindent (2) We now assume that $G \not \cong Z$. If $k \geq 2$, then putting so in the proof of Theorem~\ref{thm_lb}, we get $|V(G)| \leq 3|C|-3$ and hence, the result follows. Therefore, for the rest of this proof, let us assume that $k=1$. In other words, $G[C]$ is a connected subgraph of $G$. Now, we must have $|\mathcal{K}_{leaf}(G[C])| \geq 2$, or else, if $|\mathcal{K}_{leaf}(G[C])| = 1$, the same analysis as in (1) would follow and we would have $G \cong Z$, a contradiction. Thus, putting $|\mathcal{K}_{leaf}(G[C])| \geq 2$ in the proof of Lemma \ref{claimV4}(2), we get $|V_4| \leq |\mathcal{K}(G[C])| \leq |C|- 3k +n_1(G[C])$. Feeding this further in the proof of Theorem \ref{thm_lb} when $C$ is an OLD-code, we get $|V(G)| \leq 3|C|-3$ and hence, the lower bound holds.

\medskip
To show that there exist arbitrarily large and infinitely many block graphs whose OLD-numbers attain this bound, for $\ell \geq 6$, we consider a path on vertices $u_1, u_2, \ldots,u_\ell$. For $i=1, \ell$ and $3 \leq i \leq \ell-2$, attach a vertex $v_i$ by the edge $v_iu_i$. Moreover, for each pair $u_i,u_{i+1}$ for $1 \leq i \leq \ell-1$, attach a vertex $w_i$ by edges $w_iu_i$ and $w_iu_{i+1}$. We call the graph $G$. See Figure \ref{fig_OLD2_extremal} for an example with $\ell=6$. Then it can be verified that $C = \{u_1, u_2, \ldots , u_\ell\}$ is a minimum OLD-code of $G$. Since $|V(G)| = 3\ell-3$, the OLD-number of $G$ attains the lower bound.
\end{proof}

If we now consider the parameter $|\mathcal{K}(G)|$, we can use the relation $|V(G)|\geq |\mathcal{K}(G)|+1$ to obtain a similar lower bound. But this lower bound can be improved as the next result shows.

\begin{theorem} \label{thm_lb_K(G)}
  Let $G$ be a connected block graph and $\mathcal{K}(G)$ be the set of all blocks of $G$. Then we have
  \begin{itemize}
  \item $\gamma^{ID}(G) \geq \frac{3(|\mathcal{K}(G)|+2)}{7}$;

  \item $\gamma^{LD}(G)\geq \frac{|\mathcal{K}(G)|+2}{3}$; and

  \item $\gamma^{OLD}(G)\geq \frac{|\mathcal{K}(G)|}{2}+1$.
  \end{itemize}
  \end{theorem}

\begin{proof}
Assume $C$ to be a code of $G$ and $G[C_1], G[C_2], \ldots, G[C_k]$ be all the connected components of $G[C]$, where $C_1, C_2, \ldots , C_k$ are subsets of $C$. First, we define $\mathcal{I}_G(C) = \{ K \in \mathcal{K}(G) : V(L) \subset V(K) \text{ for some } L \in \mathcal{K}(G[C]) \}$. Moreover, for each $1 \leq i \leq k$, let $\mathcal{I}_G(C_i) = \{ K \in \mathcal{K}(G) : V(L) \subset V(K) \text{ for some } L \in \mathcal{K}(G[C_i]) \}$. Next, we define the following types of blocks of $G$.
\begin{enumerate}
  	\item Let $\mathcal{K}_C(G) = \{ K \in \mathcal{I}_G(C) : V(K) \subset C \}$, i.e. all blocks of $G$ which are also blocks of the subgraph $G[C]$ (also a block graph) of $G$.
  	
	\item Let $\mathcal{K}_{\overline{C}}(G) = \mathcal{K}(G) \setminus \mathcal{I}_G(C)$. In other words, the set $\mathcal{K}_{\overline{C}}(G)$ includes all blocks of $G$ which do not contain any vertices of the code $C$.
	
  	\item For $i=2,3,4$, let $\mathcal{K}_i(G) = \{ K \in \mathcal{I}_G(C) : V(K) \cap V_i \neq \emptyset \}$ (recall the sets $V_1, V_2, V_3, V_4$ from Definition \ref{def_V_i}).
  \end{enumerate}

We note here that, $\mathcal{K}(G) = \mathcal{K}_C(G) \cup \mathcal{K}_{\overline{C}}(G) \cup \mathcal{K}_2(G) \cup \mathcal{K}_3(G) \cup \mathcal{K}_4(G)$. We now have the following bounds.

\begin{claim}
$|\mathcal{K}_2 (G)| \leq |V_2|$.
\end{claim}

\begin{proofofclaim}
Since each vertex in the set $V_2$ belongs to a unique block $K \in \mathcal{K}_2(G)$, this claim is true.
\end{proofofclaim}

We now invoke the auxiliary graph $F_C(G)$ of $G$ from Definition \ref{def_aux} and assume that there are $l$ connected components of $F_C(G)$. Then we have the following claim.

\begin{claim}
$|\mathcal{K}_3 (G)| \leq 2(k-l)$.
\end{claim}

\begin{proofofclaim}
Since each vertex of $F_C(G)$ in the part $A$ is of degree at least $2$, we have $|E(F_C(G))|$ $\geq 2|A| = 2|V_3|$. Combining this with the fact that $|E(F_C(G))| = |V_3|+k-l$ (since $F_C(G)$ is a forest by Lemma \ref{lem_F_C(G)}), we have $|V_3| \leq k-l$. Hence, we have $|\mathcal{K}_3 (G)| \leq |E(F_C(G))| \leq 2(k-l)$.
\end{proofofclaim}

\begin{claim}
$|\mathcal{K}_{\overline{C}}(G)| \leq l-1$.
\end{claim}

\begin{proofofclaim}
Let $F_1, F_2, \ldots , F_l$ be the $l$ connected components of the auxiliary (bipartite and forest) graph $F_C(G)$ of $G$. To count $|\mathcal{K}_{\overline{C}} (G)|$, we first observe that any $K \in \mathcal{K}_{\overline{C}} (G)$, does not contain any non-articulation vertex of $K$. This is because, as $V(K) \cap C = \emptyset$, any non-articulation vertex of $K$ will remain non-dominated by the code $C$, a contradiction. This therefore implies that $K$ cannot be a leaf block of $G$ (and, in particular, the root-block of $G$ as well). Hence, every block $K \in \mathcal{K}_{\overline{C}}(G)$ has a positive articulation vertex $v_K$, say. Then we have a block $K' \in \mathcal{K}(G)$ with its layer $f(K') = f(K)+1$ such that $V(K) \cap V(K') = \{ v_K \}$. Moreover, since $V(K) \cap C = \emptyset$, for the code $C$ to dominate $v_K$, we may assume, without loss of generality, the block $K'$ to be such that there exists a vertex $v_{K'}$ of $G$ in $V(K') \cap C$. For every $K \in \mathcal{K}_{\overline{C}}(G)$, therefore, we fix such a triple $(v_K, K', v_{K'})$ for the rest of this proof.
Assume $G[C_i]$ to be the component of $G[C]$ such that $v_{K'} \in C_i$. Moreover, let $u_i$ (the vertex of the part $B$ of $F_C(G)$ corresponding to $G[C_i]$) belong to the component $F_j$, for some $1 \leq j \leq l$. Then, we associate $F_j$ with the block $K \in \mathcal{K}_{\overline{C}}(G)$. More precisely, we define the following mapping.
\begin{flalign*}
g: \mathcal{K}_{\overline{C}}(G) &\to \{F_1, F_2, \ldots , F_l\}\\
K &\mapsto F_j \text{, where $(v_K, K', v_{K'})$ is fixed, $v_{K'} \in C_i$ and $u_i \in F_j$}
\end{flalign*}
Since the vertex $v_{K'} \in C$ can belong to exactly one component $G[C_i]$ of $G[C]$ and, similarly, the vertex $u_i$ of $F_C(G)$ can belong to exactly one of its components $F_j$, the mapping $g$ is therefore well-defined. We now claim that $g$ is one-to-one. Indeed, consider any $K\in \mathcal{K}_{\overline{C}}(G)$ with $(v_K,K',v_{K'})$ associated to it, and such that $g(K)=F_j$. Assume by contradiction that there is $L\in \mathcal{K}_{\overline{C}}(G)$, $K\neq L$, with $g(L)=F_j$, and $(v_L,L',v_{L'})$ is associated with $L$.
Since $V(K) \cap C = \emptyset$, for every vertex $u_r\in V(F_j)$ and any block $J \in \mathcal{I}_G(C_r)$, we have $f(J)\geq f(K)+1$. Thus, $f(K')=f(K)+1$ is minimum among all such blocks. By the same reasoning applied to $L$, we also have $f(L')=f(L)+1$ minimum among all blocks $J \in \mathcal{I}_G(C_r)$, where $u_r \in V(F_j)$. This implies that $f(L')=f(K')$, since both values are minimim among all $f(J)$, where $J \in \mathcal{I}_G(C_r)$ and $u_r \in V(F_j)$. Therefore, we also have $f(K) = f(L)$. Now, since $K \neq L$, we must have $V(K) \cap V(L) = \emptyset$, or else, we would have $|f(K) - f(L)|=1$, a contradiction. Now, let $G_K$ represent the \emph{descendant block graph} of $K$, that is, the connected subgraph (also a block graph) of $G$ rooted at $K$. Similarly define $G_L$ to be the descendant block graph of $L$. Then, by the structure of the block graph $G$, its subgraphs $G_K$ and $G_L$ are vertex-disjoint. Since $V(K) \cap C = V(L) \cap C = \emptyset$, all components $G[C_r]$ of $G[C]$ for all $u_r \in V(F_j)$ must belong to $G_K$ and $G_L$ simultaneously, which is a contradiction since $V(G_K) \cap V(G_L) = \emptyset$.
This shows that $g$ is one-to-one, as claimed.

\medskip
As noticed before, for each $K \in \mathcal{K}_{\overline{C}}(G)$ such that $g(K) = F_j$, we have
$$f(K)+1 = \min \{f(J) : J \in \mathcal{I}_G(C_r), u_r \in F_j\}.$$
Now, let $K_0$ be the root block of $G$. As argued before, $K_0 \notin \mathcal{K}_{\overline{C}}(G)$. In other words, $V(K_0) \cap C_{i_0} \neq \emptyset$ for some component $G[C_{i_0}]$ of $G[C]$. Let $u_{i_0} \in F_{j_0}$. Therefore, we have
$$0 \leq \min \{f(J) : J \in \mathcal{I}_G(C_r), u_r \in F_{j_0} \} \leq f(K_0) = 0 \quad \text{(since $K_0 \in \mathcal{I}_G(C_{i_0})$)}.$$
This implies that $g(K) = F_{i_0}$ for any $K$ would imply $f(K) = -1$, which is not possible. Hence, the image of the function $g$ is a subset of $\{F_1, F_2, \ldots , F_l\} \setminus \{F_{i_0}\}$. This, along with the fact that $g$ is one-to-one, therefore implies that $|\mathcal{K}_{\overline{C}}(G)| \leq l-1$.
\end{proofofclaim}

\begin{claim}
$|\mathcal{K}_C(G) \cup \mathcal{K}_4(G)| \leq |\mathcal{K}(G[C])| - n_0(G[C])$.
\end{claim}

\begin{proofofclaim}
Assume that $K \in \mathcal{K}(G)$ is a block of $\mathcal{K}_C(G) \cup \mathcal{K}_4(G)$. Then, $V(K)$ contains at least two vertices, say, $u, v \in C$. Therefore, $uv \in E(G)$. So, assume $L \in \mathcal{K}(G[C])$ to be the~block such that $u, v \in V(L)$. Then, $V(L) \subset V(K)$. Thus, every block $K \in \mathcal{K}_C(G) \cup \mathcal{K}_4(G)$ can be associated with a block $L \in \mathcal{K}(G[C])$ such that $|V(L)| \geq 2$. Moreover, by the structure of a block~graph, this association is one-to-one. This implies that \mbox{$|\mathcal{K}_C(G) \cup \mathcal{K}_4(G)| \leq |\mathcal{K}(G[C])| - n_0(G[C])$.}
\end{proofofclaim}

To compute $|\mathcal{K}(G)|$ now, we have from the above Claims A, B, C and D that
\begin{flalign} \label{eqn K(G)}
|\mathcal{K}(G)| &\leq |\mathcal{K}_C(G) \cup \mathcal{K}_4(G)| + |\mathcal{K}_{\overline{C}}(G)| + |\mathcal{K}_2(G)| + |\mathcal{K}_3(G)| \nonumber\\
&\leq |\mathcal{K}(G[C])| - n_0(G[C]) + l-1 + |V_2| + 2(k-l).
\end{flalign}
Therefore, using Equation (\ref{eqn K(G)}), we have the following.

\medskip
For ID-codes:
\begin{flalign*}
|\mathcal{K}(G)| &\leq |\mathcal{K}(G[C])|- n_0(G[C]) + l-1 + |V_2| + 2(k-l)\\
&\leq |C|-k + l-1 + |C|-n_0(G[C]) + 2(k-l) &&\hspace*{-16mm}\text{[using Corollary \ref{cor_K(G) <= V-k} and Lemma \ref{claimV2}(1)]}~~\\
&= 2|C| + k-l - n_0(G[C]) - 1\\
&\leq 2|C| + k - n_0(G[C]) - 2.
\end{flalign*}

Now, $k-n_0(G[C])$ is the total number of components of $G[C]$ of order at least $2$. Any such component must contain at least $3$ vertices of the code $C$ (or else, if a component contained exactly two vertices, they would not be closed-separated in $G$ by $C$). Therefore, $3(k-n_0(G[C])) \leq |C| - n_0(G[C])$. Therefore, $3|\mathcal{K}(G)| \leq 7|C|-n_0(G[C])-6 \leq 7|C|-6$ and, hence, the result holds.

\medskip
For LD-codes:
\begin{flalign*}
|\mathcal{K}(G)| &\leq |\mathcal{K}(G[C])| - n_0(G[C]) + l-1 + |V_2| + 2(k-l)\\
&\leq |C|-k + l-1 + |C| + 2(k-l) && \hspace*{-16mm}\text{[using Corollary \ref{cor_K(G) <= V-k} and Lemma \ref{claimV2}(2)]}~~\\
&= 2|C| + k-l - 1\\
&\leq 3|C| - 2.
\end{flalign*}

Finally, for OLD-codes, we have $n_0(G[C])=0$. Hence,
\begin{flalign*}
|\mathcal{K}(G)| &\leq |\mathcal{K}(G[C])| + l-1 + |V_2| + 2(k-l)\\
&\leq |C|-2k + l-1 + |C| + 2(k-l) &&\hspace*{-10mm}\text{[using Lemmas \ref{claimV2}(3) and \ref{claimV4}(2)]}~~\\
&=2|C| - l-1 \\
& \leq 2|C| - 2, \quad \text{using $l \geq 1$}.
\end{flalign*}
This proves the theorem.
\end{proof}

Recall that a block in a graph $G$ (not necessarily a block graph) is a maximal complete subgraph of $G$. Then, the lower bounds in Theorem~\ref{thm_lb_K(G)} in terms of the number of blocks in a graph fail for graphs in general. For example, consider the split graph $G$ with its vertex set $V(G)=\{v_1,...,v_k\}\cup\{u_X : X \subseteq \{1,...,k\} \text{ and } X\neq \emptyset\}$. The vertices $v_1,\ldots, v_k$ induce a clique, whereas the vertices $u_X$ induce an independent set. Moreover, there is an edge between the vertices $v_i$ and $u_X$ if and only if $i\in X$. This graph has an identifying code of size $2k$ (the clique with the vertices corresponding to the singletons), but the number of blocks in $G$ is $2^k$. Similar examples exist to show that the LD- and OLD-numbers of some split graphs that are not block-graphs violate the lower bounds in Theorem~\ref{thm_lb_K(G)}.

\medskip
Note that, for any tree $G$, we have $|\mathcal{K}(G)|=|E(G)|=|V(G)|-1$. Thus, for trees (which are particular block graphs with each block being of order $2$), Theorem \ref{thm_lb_K(G)} provides the same lower bound \Big ($\frac{3(|V(G)|+1)}{7}$\Big) for ID-numbers as was given in \cite{bertrand20051}. The lower bounds for LD-numbers of trees given in Theorems \ref{thm_lb} and \ref{thm_lb_K(G)} are the same; and which, in turn, are the same as that given in \cite{slater1987domination}. The lower bound for OLD-numbers of trees given by Theorem \ref{thm_lb_K(G)} is the same when $|V(G)|$ is even and one short when $|V(G)|$ is odd as the lower bound given in \cite{seo2010open} \Big (which is $\Big \lceil \frac{|V(G)|}{2} \Big \rceil +1$\Big). However, for $|V(G)|$ odd, this lower bound for OLD-numbers given by Theorem \ref{thm_lb_K(G)} is tight for block graphs in general (see Figure \ref{fig_OLD_extremal} and E.g. (3) below). The following are examples of extremal block graphs described in Propositions \ref{prop_ID_extremal}, \ref{prop_LD_extremal} and \ref{prop_OLD_extremal} whose code numbers attain the lower bounds in Theorem~\ref{thm_lb_K(G)}.
\begin{enumerate}[E.g. (1)]
\item For $k=3$ in Proposition \ref{prop_ID_extremal}, the graph $G$ is a $1$-corona of a $P_3$ (with $|\mathcal{K}(G)|=5$) and is an extremal example whose ID-number ($=3$) attains the bound in Theorem \ref{thm_lb_K(G)}.

\item In the proof of Proposition \ref{prop_LD_extremal}, the graph $G'$ obtained by deleting all the edges $u_i u_{i+1}$ of $G$ is an example whose LD-number ($=k$) attains the bound in Theorem \ref{thm_lb_K(G)} (note that $|\mathcal{K}(G')|=3k-2$).

\item Finally, the graph $Z$ described in Proposition \ref{prop_OLD_extremal}(1) (with $|\mathcal{K}(G)|=4$) also serves as an example whose OLD-number ($=3$) and  attains the bound given in Theorem \ref{thm_lb_K(G)}.
\end{enumerate}

Apart from the above examples, there are infinite families of trees reaching the three bounds in Theorem \ref{thm_lb_K(G)} (see \cite{bertrand20051} for ID-codes, \cite{slater1987domination} for LD-codes and \cite{seo2010open} for OLD-codes).

\section{Conclusion}\label{sec:conclu}

Block graphs form a subclass of chordal graphs for which all three considered identification problems can be solved in linear time~\cite{argiroffo2020linear}. In this paper, we complemented this result by presenting lower and upper bounds for all three codes.
We gave bounds using both the number of vertices — as it has been done for several other classes of graphs — and also using the parameter $|\mathcal{K}(G)|$ of the number of blocks of $G$ that is more fitting for block graphs. In particular, we verified a conjecture from~\cite{ABLW_ICGT} (Conjecture~\ref{Conj_ID}) concerning an upper bound for $\gamma^{ID}(G)$, and also proved the conjecture on the LD-number~\cite{garijo2014difference} (Conjecture~\ref{Conj_LD twin free}) for the special case of block graphs. Moreover, we addressed the questions to find block graphs where the provided lower and upper bounds are attained.

\medskip
The structural properties of block graphs have enabled us to prove interesting bounds for the three considered identification problems. It would be interesting to see whether other structured classes can be studied in a similar way. It would also be interesting to prove Conjecture~\ref{Conj_LD twin free} for a larger class of graphs, for example for all chordal graphs.

\subsection*{Acknowledgments}

\noindent Works by Dipayan Chakraborty and Florent Foucaud were partially supported by the French government initiative IDEX-ISITE CAP 20-25 (ANR-16-IDEX-0001), the International Research Center ``Innovative Transport and Production Systems" I-SITE CAP 20-25 and the
ANR GRALMECO project (ANR -21- CE48-0004).



\end{document}